\documentclass{article}
\usepackage{authblk}
\usepackage{url}

\usepackage{amsmath,amsfonts}
\usepackage{halloweenmath}[2019/11/01]
\newcommand{\pad}[2]{\frac{\partial #1}{\partial #2}}
\newcommand{\RR}{\mathbb{R}}
\usepackage[percent]{overpic}
\newcommand{\nodes}{{\mathcal A}}
\newcommand{\vphi}{\varphi}

\usepackage{todonotes}
\usepackage{float} 
\usepackage{booktabs}

\begin{document}

%

\title{A comparison of the Coco-Russo scheme \\ and $\protect\mathghost$-FEM for elliptic equations \\ in arbitrary domains}
%
%
\author[1,2]{Clarissa Astuto} 
\author[2]{Armando Coco}
\author[3]{Umberto Zerbinati} 
%
%
%
\affil[1]{Department of Applied Mathematics and Computational Sciences, King Abdullah University of Science and Technology (KAUST), Saudi Arabia}
\affil[2]{Department of Mathematics and Computer Science, University of Catania, Italy}
\affil[3]{Mathematical Institute, University of Oxford, United Kingdom}

\maketitle              

\begin{abstract}
In this paper, a comparative study between the Coco-Russo scheme (based on finite-difference scheme) and the $\mathghost$-FEM (based on finite-element method) is presented when solving the Poisson equation in arbitrary domains. The comparison between the two numerical methods is carried out by presenting analytical results from the literature \cite{cocoStissi,astuto2024nodal}, together with numerical tests in various geometries and boundary conditions.
\end{abstract}

\section{Introduction}
Elliptic partial differential equations (PDEs) are used to describe a large variety of physical phenomena.
In this paper, we will focus on the prototypical elliptic PDE, i.e. the Poisson equation with mixed boundary conditions,
\begin{equation}
	\label{eq:poissonMixed}
	-\Delta u = f \; {\rm in }\;\Omega, \; u=g_D\;{\rm on }\;\Gamma_D \subset \partial\Omega \; {\rm and } \;\pad{u}{\widehat n}  = g_N \; {\rm on }\;\Gamma_N\subset\partial\Omega,
\end{equation}
where $\Omega\subset R\subseteq\RR^2$, $R$ is a rectangular region, $\Gamma_D \cup \Gamma_N = \partial \Omega =:\Gamma,$ with $\Gamma_D \cap \Gamma_N = \emptyset$ and $\widehat n$ is the outgoing normal vector to $\Gamma$. 

Many numerical schemes have been developed to solve the Poisson equation with mixed boundary conditions on an arbitrary domain $\Omega$.
Some of these numerical schemes will be based on the strong formulation of the problem, while others will be based on reformulating \eqref{eq:poissonMixed} as a variational problem.
Primary examples of the numerical methods based on the strong formulation are unfitted finite difference schemes, which rely on the description of the domain by a level set function  \cite{Gibou2002,Gibou2007,Shortley1938}. Among such methods, we will focus on the Coco-Russo scheme \cite{Coco2013,Coco2018} which has been applied in several contexts~\cite{COCO2020109623,astuto2023multiscale,ASTUTO2023111880,astuto2023time,astuto2024high}.
On the other hand, Galerkin methods are based on the variational formulation of the problem, and they have the advantage of preserving more structural properties of the continuous problem. A prime example of Galerkin method is the finite element method (FEM), which relies on a reformulation of \eqref{eq:poissonMixed} as a variational problem and on discrete space constructed from a ``good'' tessellation of the domain, \cite{BrennerScott,Ciarlet}.
An interesting variant of the finite element method is the cut finite element method (CutFEM), where a level set function is used to represent the boundary of $\Omega$, \cite{Burman2015,Lehrenfeld2016,Lehrenfeld2018,Burman2022}.
In this work, we will focus our attention on the $\mathghost$-FEM, a ghost nodal method based on variational formulation, which is a variant of the finite-element method that uses level set functions to represent the boundary of $\Omega$, similar to the CutFEM, \cite{astuto2024nodal}.

Although the unfitted finite difference method may appear simpler to implement and understand for novice users, finite element methods not only have the advantage of preserving more structural properties of the continuous problem but also the theoretical framework behind such schemes allows to prove more complete convergence results. 
This paper aims to provide a comparative study between the Coco-Russo scheme and the $\mathghost$-FEM, when solving the Poisson equation with mixed boundary conditions.

{
The paper is structured as follows. In Section~\ref{sect:numschemes}, we describe the two numerical schemes of interest: Coco-Russo scheme and $\mathghost$-FEM. Section~\ref{sec:remarks} first discusses \textit{a priori} error estimates for both numerical methods. It then addresses the conditioning of these methods in critical situations, such as the presence of small cells. Finally, it presents a subsection on fast solvers. In Section~\ref{sect:results}, we test the two approaches in several domains, providing a discussion on their respective strengths and weaknesses. At the end we draw some conclusions.
}

\section{Numerical schemes}\label{sect:numschemes}
\label{sec:numerical_schemes}
The present section will be devoted to the description of the numerical schemes used in this work, namely the Coco-Russo scheme and the $\mathghost$-FEM.
Let the domain $\Omega \subset R$ be described by a level set function $\phi(x,y)$ that is positive inside $\Omega$, negative in $R\setminus \Omega$ and zero on the boundary $\Gamma$ (see, for example, \cite{sussman1994level,Osher,russo2000remark,book:72748}):
\begin{eqnarray}
	\Omega = \{(x,y): \phi(x,y) > 0\}, \qquad
	\Gamma = \{(x,y): \phi(x,y) = 0\}.
\end{eqnarray}
The unit normal vector $\widehat n$ in \eqref{eq:poissonMixed} can be computed as $\widehat n= \frac{\nabla \phi }{|\nabla \phi|}$
where the level-set function $\phi$ is assumed to be explicitly known.

\begin{figure}[htp]
	\centering
	\begin{minipage}
		{.32\textwidth}	\centering
\begin{overpic}[abs,width=\textwidth,unit=1mm,scale=.25]{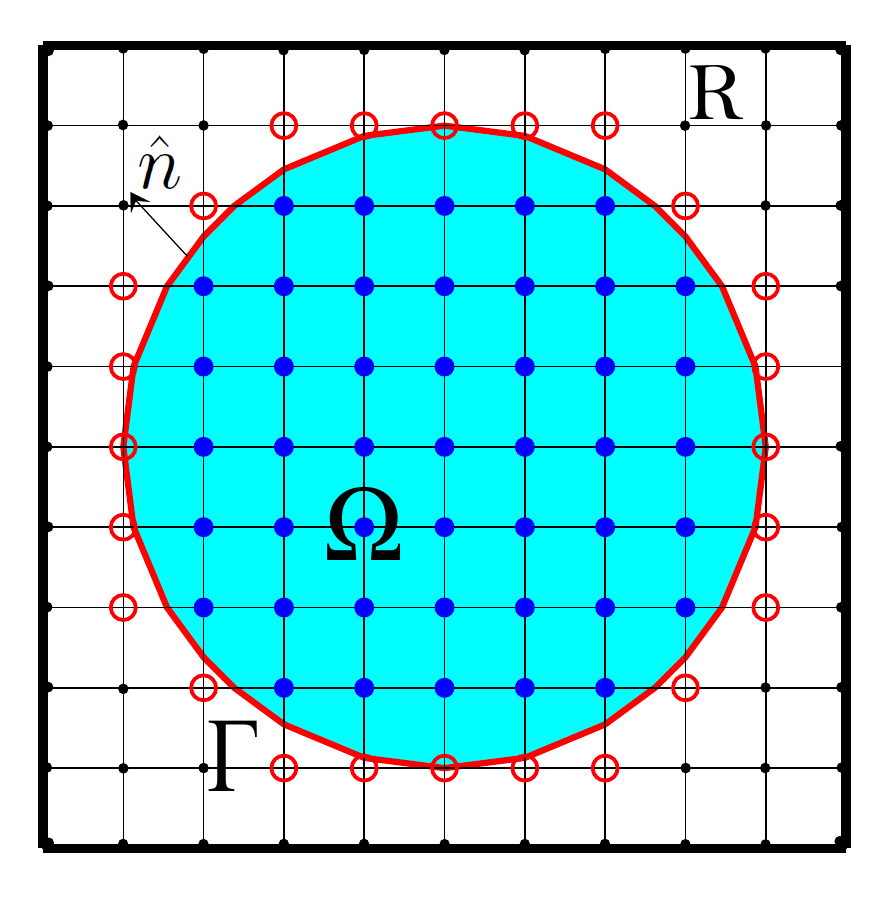}
\put(-4.,37){(a)}
\end{overpic}
	\end{minipage}
	\begin{minipage}
		{.32\textwidth}	\centering
\begin{overpic}[abs,width=\textwidth,unit=1mm,scale=.25]{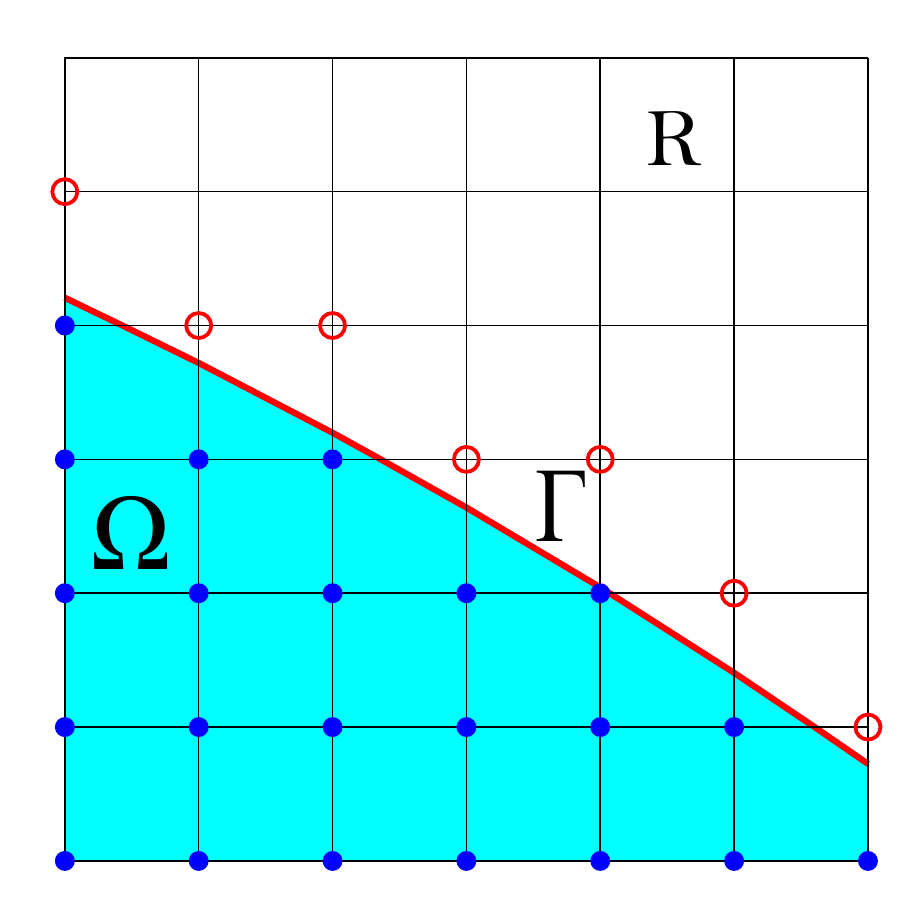}
\put(-2.,37){(b)} 
\put(12.,38){Coco-Russo} 
\end{overpic}
	\end{minipage}
	\begin{minipage}{.32\textwidth}	
  \centering
\begin{overpic}[abs,width=\textwidth,unit=1mm,scale=.25]{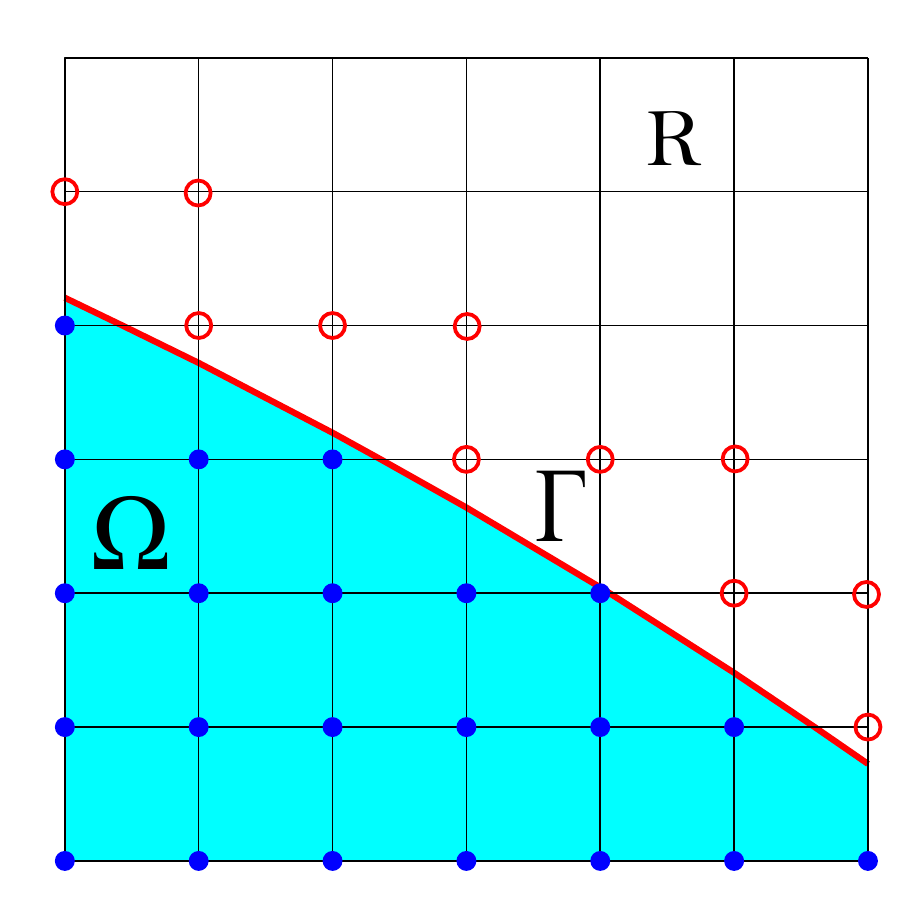}
\put(-2.,37){(c)} 
\put(14.,38){$\protect \mathghost$-FEM} 
\end{overpic}
	\end{minipage}
	\caption{\textit{(a): Representation of the domain $\Omega\subset R$ and of the normal vector $\widehat n$ to the boundary $\Gamma$. The classification is the following: internal points (blue points),  ghost points (red circles), and inactive points (small black dots). Furthermore, we show the different distribution of ghost points for the Coco-Russo scheme in panel (b) and for the $\protect\mathghost$-FEM in (c).}}
	\label{fig:levelset}
\end{figure}

\subsection{Coco-Russo scheme}
Let us consider a uniform square Cartesian cell-centered discretization, with $\Delta x = \Delta y =: h$, and the set of grid points is $\mathcal{S}_h = (x_h,y_h) = \{(x_i,y_j)=(ih,jh), (i,j) \in \{0,\cdots,N\}^2 \}$, where $N \in \mathbb{N}$ and 
$h = L_x/N$ with $L_x = L_y = 2$. Here we define the set of internal points $\Omega_h = \mathcal{S}_h \cap \Omega$, 
and the set of ghost points $\mathcal{G}_h$ {as the external points with at least one neighbor internal point}, formally defined as follows
\begin{equation}
	(x_i,y_j) \in \mathcal{G}_h \iff (x_i,y_j) \in \mathcal{S}_h \text{ and } \{(x_i \pm h,y_j),(x_i,y_j\pm h) \} \cap \Omega_h \neq \emptyset.
\end{equation}
The other grid points, $\mathcal{S}_h\setminus(\Omega_h \cup \mathcal{G}_h)$, are called inactive points. See Fig.~\ref{fig:levelset} (a) for a classification of inside, ghost, and inactive points.
Let $N_I = |\Omega_h|$ and $N_G = |\mathcal{G}_h|$ be the cardinality of the sets $\Omega_h$ and $\mathcal{G}_h$, respectively, and $\mathcal{N} = N_I + N_G$ the total number of active points. 
We compute the solution $ u_h $ at the grid points of $\Omega_h \cup \mathcal{G}_h$, using a finite difference discretization of the equations on the $N_I$ internal grid points together with a suitable interpolation
to define the $N_G$ ghost values from the boundary conditions. Since each equation on a ghost point may involve other ghost points, the equations on the ghost points might be coupled. For this reason, the whole $\mathcal{N} \times \mathcal{N}$ system, with non-eliminated boundary conditions, is considered. 
\begin{figure}[htp]
	\centering
\begin{overpic}[abs,width=0.6\textwidth,unit=1mm,scale=.25]{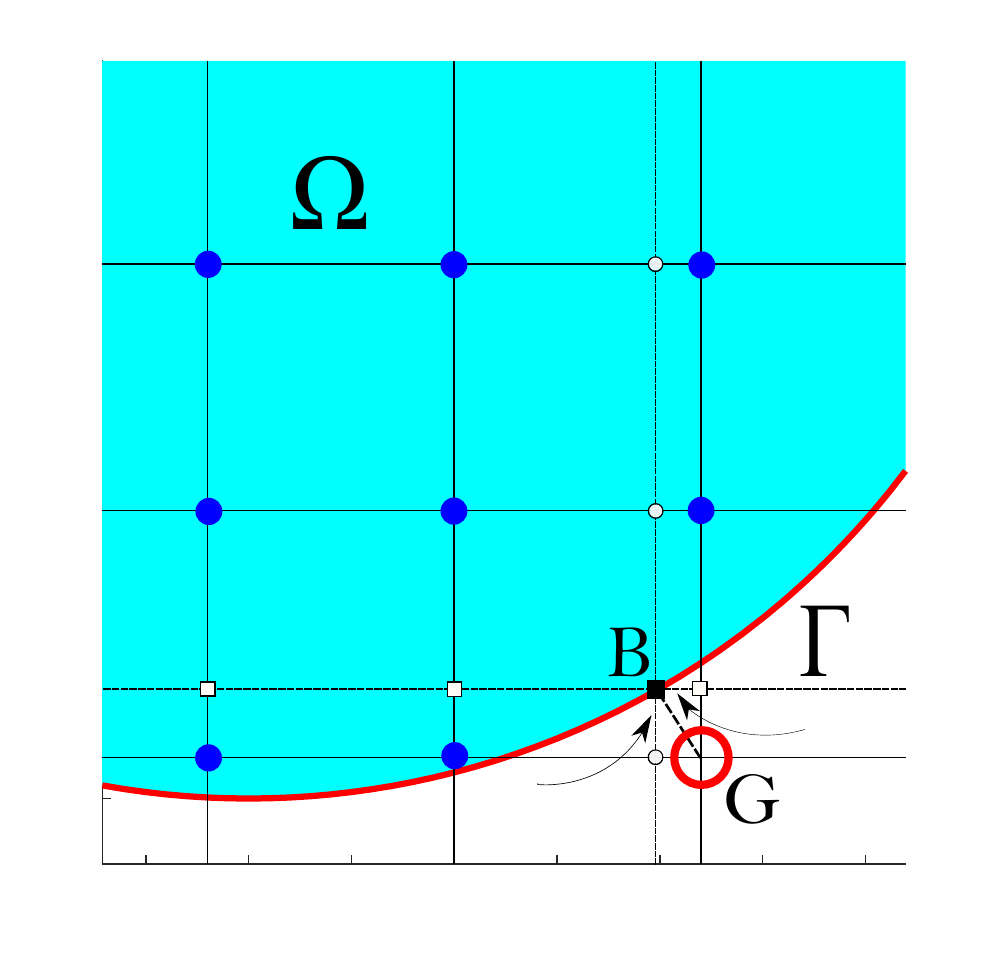}
\put(36,9){\small $\vartheta_y h$}
\put(59,16){\small $\vartheta_x h$}
\end{overpic}
\caption{\textit{We show the nine-point stencil for the interpolation operator for Coco-Russo scheme. G (red circle) is the ghost point, B (black square) is the closest point to G that belongs to $\Gamma$ and the eight blue points complete the nine-point stencil. The empty black circles are the points of interpolation in x direction and the empty black squares the ones in y direction. The two quantities $\vartheta_x$ and $\vartheta_y$ are defined in Eq.\eqref{expr_theta}.}}
	\label{fig:stencil}
\end{figure}

{The discretization of the problem \eqref{eq:poissonMixed} leads to a linear system}
\begin{eqnarray}
\label{eq:LS_FDM}
	A^{\rm FDM}  u_{h} = f_h,
\end{eqnarray}
where $A^{\rm FDM}$ is a $\mathcal{N} \times \mathcal{N}$ matrix representing the discretization of the derivatives and the interpolation operators. 
If $P_{ij} = (x_i,y_j) \in \Omega_h$ is an internal grid point (as in Fig.~\ref{fig:levelset} (a)), we discretize the {elliptic operator by the standard central finite-difference formula}
\begin{eqnarray} 
\label{eq_Lh}
A^{\rm FDM}\, u_h\Big|_{i,j} &=&  \frac{u_{i,j+1}  + u_{i,j-1} + u_{i+1,j}  + u_{i-1,j} - 4u_{i,j}}{h^2}.
\end{eqnarray}
If $G=(x_i,y_j) \in \mathcal{G}_h$ is a ghost point, then we discretize the boundary condition on $\Gamma$, following the ghost-point approach proposed in \cite{Coco2013,Coco2018}, and summarised as follows. 
Initially, we determine the closest boundary point to $G$, denoted as $B (\in \Gamma)$, utilizing the formula:
\[
B = G - \widehat{n}_G \cdot \nu,
\]
where $\widehat{n}_G$ is computed by $ \widehat{n}_G = \nabla \phi / |\nabla \phi| $ and $ \nabla \phi $ is discretised via standard finite differences centered at $G$. The parameter $\nu$  is determined by solving the equation $\phi(G - \widehat{n}_G \cdot \nu) = 0 $, via a bisection method with a tolerance of \( 10^{-4} h \) for the distance between the \( B \) and the boundary \( \Gamma \). Within the bisection method, the evaluation of \( \phi \) at off-grid points is accomplished through bilinear interpolation.

Then, we identify the upwind nine-point stencil starting from $G=(x_G,y_G)=(x_i,y_j)$, containing $B=(x_B,y_B)$:
\[
\left\{ (x_{i+s_x m_x},x_{j+s_y m_y}) \colon m_x,m_y=0,1,2 \right\},
\]
where $s_x = \frac{x_B-x_G}{\lvert{x_B-x_G}\lvert}$ and $s_y = \frac{y_B-y_G}{\lvert y_B-y_G\lvert}$. The solution $ u_h $ and its first derivative are then interpolated at the boundary point $B$ using the discrete values $ u_{i,j}$ on the nine-point stencil.
In particular, we start defining (see Fig.~\ref{fig:stencil})
\begin{align}
\label{expr_theta}
\vartheta_x =  s_x (x_B-x_G)/h, \qquad
\vartheta_y =  s_y (y_B-y_G)/h,
\end{align}
with $0\leq \vartheta_x,\vartheta_y < 1$.
For a generic grid function $c$, its interpolant $\widetilde{c}$ and its partial derivatives are evaluated  at $B$ by:

\begin{eqnarray}
    \label{coeffsLSstencil} 
\widetilde{c}(B) &= \sum_{m_x,m_y=0}^p l_{m_x}(\vartheta_x) l_{m_y}(\vartheta_y)  c_{i+s_x m_x,j+s_y m_y},
\\ \label{coeffsLSstencil2}
\frac{\partial \widetilde{c}}{\partial x}(B) &= s_x \sum_{m_x,m_y=0}^p l'_{m_x}(\vartheta_x) l_{m_y}(\vartheta_y)  c_{i+s_x m_x,j+s_y m_y},
\\ \label{coeffsLSstencil3}
\frac{\partial \widetilde{c}}{\partial y}(B) &= s_y \sum_{m_x,m_y=0}^2 l_{m_x}(\vartheta_x) l'_{m_y}(\vartheta_y)  c _{i+s_x m_x,j+s_y m_y},
\end{eqnarray}
The normal derivative at $B$ is computed via
\begin{equation}\label{eq:normal_deriv}
\frac{\partial \widetilde{c}}{\partial n}(B) = 
\nabla \widetilde{c} (B) \cdot \widehat{n}_B, \qquad \text{ with }
\widehat{n}_B = \frac{\nabla \widetilde{\phi} (B) }{ \left| \nabla \widetilde{\phi} (B) \right|}.
\end{equation}
The integer $p$ represents the size of the interpolation stencil.
For $p = 1$ we obtain a 4-point stencil, and the coefficients are
{
\begin{equation}\label{coeffs_l1}
l(\vartheta_\mu)=\left(1-\vartheta_\mu, \,\vartheta_\mu\right),
\qquad
l'(\vartheta_\mu) = \frac{1}{h}\left( -1, \, 1\right),  \quad \mu = x,y.
\end{equation}
}
{
With this choice, the method achieves second-order accuracy for Dirichlet boundary conditions but only first-order accuracy for Neumann boundary conditions~\cite{Coco2013}. Additionally, the gradient of the solution is only first-order accurate for both Dirichlet and Neumann boundary conditions. To achieve second-order accuracy for both the solution and its gradient under Dirichlet and Neumann boundary conditions,
we must increase the size of the interpolation stencil from 4 to 9 points, namely $p=2$, and the coefficients are
}
\begin{align}
\label{coeffs_l2}
l(\vartheta_\mu) &\!=\! \left( \frac{(1\!-\vartheta_\mu)(2\!-\vartheta_\mu)}{2}, \, \vartheta_\mu (2\!-\vartheta_\mu), \, \frac{\vartheta_\mu(\vartheta_\mu\!-\!1)}{2} \right),
\\
l'(\vartheta_\mu) &= \frac{1}{h} \left( \frac{(2\vartheta_\mu-3)}{2}, \, 2(1-\vartheta_\mu), \, \frac{(2\vartheta_\mu-1)}{2} \right),  \; \mu \!=\! x,y.
\end{align}

{
We observe that the method can be extended to high-order by increasing the stencil size. In~\cite{CocoHighOrder} the Coco-Russo method has been extended to fourth order of accuracy.
}
Finally, the rows of linear system $A^{\rm FDM} u_h = f_h$ associated to the ghost point $G=(x_G,y_G)$ are defined by evaluating the boundary condition in $B \in \Gamma$, i.e.
\begin{eqnarray}\label{QHghost}
	\left( A^{\rm FDM} u_h = f_h \right)\Big|_G \text{ is obtained from }\widetilde{u}_h(B) = g_D(B), \qquad \text{ if } B \in \Gamma_D \\ \label{QHghost2}
 	\left( A^{\rm FDM} u_h = f_h \right)\Big|_G \text{ is obtained from } \frac{\partial \widetilde{u}_h}{\partial n}(B) = g_N(B), \qquad \text{ if } B\in \Gamma_N.
\end{eqnarray}

If the domain is a circle, the computation of the normal direction is exact: $\widehat{n}_G = ({G-O})/{|G-O|},$ where $O$ is the center of the circle. We use the exact formula for the numerical tests with circular domains in Section~\ref{sect:results}.
\subsection{Ghost nodal Finite Element method ($\protect\mathghost$-FEM)}
We now consider the variational formulation of problem \eqref{eq:poissonMixed}, i.e. we look for $u \in H^1(\Omega)$ such that
\begin{equation}
	\int_\Omega\nabla u\cdot \nabla v\, d\Omega 
	 +  \int_{\Gamma_D} \Big(\lambda\cdot (u-g_D)  -\pad{u}{n}\Big)v \,d\Gamma
	= \int_\Omega f v \, d\Omega 
	+  \int_{\Gamma_N} g_N v \,d\Gamma,  
	\label{eq:weak2M}
\end{equation}
$\forall v \in H^1(\Omega)$, where we made use of Neumann boundary conditions by replacing $\displaystyle \frac{\partial u}{\partial n}$ with $g_N$ on $\Gamma_N$, and of a penalisation term to approximately impose Dirichlet conditions on the boundary $\Gamma_D$.
The above problem is equivalent to impose Robin boundary conditions that as $\lambda$ approaches infinity degenerate to Dirichlet boundary conditions, \cite{Courant}.

Since the variational formulation \eqref{eq:weak2M} is not symmetric, we resort Nitsche's strategy to symmetries the variational formulation, i.e. we add and subtract the term $\int_{\Gamma_D} u \pad{v}{n} \,d\Gamma$  and look for $u \in H^1(\Omega)$ such that the following holds for all $v \in H^1(\Omega)$:
\begin{align}
    \int_\Omega\nabla u\cdot \nabla v\, d\Omega 
	 +   \int_{\Gamma_D} \Big(\lambda\cdot (u-g_D) -\pad{u}{n}\Big)v \,d\Gamma - \int_{\Gamma_D}  u \pad{v}{n} \,d\Gamma \nonumber \\
	=  
	 \int_\Omega f v \, d\Omega 
	+  \int_{\Gamma_N} g_N v \,d\Gamma - \int_{\Gamma_D} g_D \pad{v}{n} \,d\Gamma.
	\label{eq:weak2M_Nitsche}
\end{align}

Notice that, for purely Neumann problems there is no penalisation term, since $\Gamma_N = \Gamma$ and $\Gamma_D = \emptyset$, and an additional  compatibility condition has to be imposed on the solution to ensure well-posedness, i.e.
\begin{equation}
    \int_\Omega f\, d\Omega + \int_\Gamma g_N\, d\Gamma = 0.
\end{equation}
The $\mathghost$-FEM discretization of the problem is obtained by replacing the domain $\Omega$ with its polygonal approximation $\Omega_h$, whose boundary is $\Gamma_h$.
The set of grid points will be denoted by $\mathcal N$, with $\# \mathcal N = (1+N)^2$, the set of active nodes (i.e.\ internal $\mathcal{I}$ or ghost $\mathcal{G}$) will be denoted by $\mathcal{I}\cup\mathcal{G} = \nodes \subset \mathcal N$, while the set of inactive points will be referred to as $\mathcal O \subset \mathcal N$, with $\mathcal O\cup\mathcal A = \mathcal N$ and $\mathcal O \cap \mathcal A = \emptyset$ and the set of cells by $\mathcal C$, with $\# \mathcal C = N^2$. {Finally, we denote by $\Omega_c = R\setminus \Omega$ the outer region in $R$.}

Here, we define the set of ghost points $\mathcal{G}$, which are {grid} points that belong to $\Omega_c$, with at least an internal point as neighbor, formally defined as
\begin{equation}
\notag
	(x,y) \in \mathcal{G} \iff (x,y) \in {\mathcal N}\cap \Omega_c  \text{ and } \{(x \pm h,y),(x,y\pm h), (x \pm h,y\pm h) \} \cap \mathcal I \neq \emptyset.
\end{equation}

We then consider the discrete space $V_h$ given by the piecewise bilinear functions which are continuous in $R$.
In particular, as basis functions, we choose the following functions:
\begin{equation}
    \varphi_{i}(x,y) = \max\left\{
        \left(1-\frac{|x-x_i|}{\Delta x}\right) 
        \left(1-\frac{|y-y_i|}{\Delta y}\right),0
    \right\},
    \label{eq:V_h2}
\end{equation}
with $i = (i_1,i_2)$ the index that identifies a node on the grid. Any discrete function $u_h\in V_h$ can be represented as
\begin{equation}
    u_h(x,y) = \sum_{i\in\nodes}u_i\vphi_i(x,y).
    \label{eq:u_h2}
\end{equation}
The discretization of the problem \eqref{eq:weak2M_Nitsche} is obtained by replacing 
$u$ and $v$ by $u_h$ and $v_h$, both in $V_h$, and evaluating the integral  over its polygonal approximation $\Omega_h$ rather than over $\Omega$, as follows
\begin{align}
    &\sum_{j\in \mathcal N} u_j(\nabla\varphi_i,\nabla\varphi_j)_{L^2(\Omega_h)}  - \sum_{j\in \mathcal N} u_j\left(\pad{\varphi_i}{n},\varphi_j\right)_{L^2(\Gamma_{D,h})}  - \sum_{j\in \mathcal N} u_j\left(\varphi_j,\pad{\varphi_i}{n}\right)_{L^2(\Gamma_{D,h})} \nonumber \\ + &\lambda\sum_{j\in \mathcal N} u_j(\varphi_j,\varphi_i)_{L^2(\Gamma_{D,h})} = \sum_{j\in \mathcal N} f_j (\varphi_j,\varphi_i)_{L^2(\Omega_h)} + \sum_{j\in \mathcal N} {g_N}_j(\varphi_j,\varphi_i)_{L^2(\Gamma_{N,h})} \label{eq:weak2M_Nitsche_discrete} \\ \nonumber - & \sum_{j\in \mathcal N} {g_D}_j\left(\varphi_j,\pad{\varphi_i}{n}\right)_{L^2(\Gamma_{D,h})} + \lambda\sum_{j\in \mathcal N} {g_D}_j(\varphi_j,\varphi_i)_{L^2(\Gamma_{D,h})}   
\end{align}
where the quantities $f_i$, $g^D_i$, $g^N_i$ denote the nodal values of the source and of the Dirichlet and Neumann boundary functions, respectively. Here we rewrite \eqref{eq:weak2M_Nitsche_discrete}, as follows
\begin{align} \nonumber
    &\sum_{j\in \mathcal N} u_j\left(\underbrace{(\nabla\varphi_j,\nabla\varphi_i)_{L^2(\Omega_h)}}_{\rm S} - \underbrace{\left(\left(\pad{\varphi_j}{n},\varphi_i\right)_{L^2(\Gamma_{D,h})}  + \left(\varphi_j,\pad{\varphi_i}{n}\right)_{L^2(\Gamma_{D,h})} \right)}_{\rm S_T} \right. \\ + & \left. \lambda\underbrace{(\varphi_j,\varphi_i)_{L^2(\Gamma_{D,h})}}_{\rm P_{\Gamma_D}} \right)  = \sum_{j\in \mathcal N} f_j \underbrace{(\varphi_j,\varphi_i)_{L^2(\Omega_h)}}_{\rm M} + \sum_{j\in \mathcal N} {g_N}_j\underbrace{(\varphi_j,\varphi_i)_{L^2(\Gamma_{N,h})} }_{\rm N_{\Gamma_N}} \\ &+ \sum_{j\in \mathcal N} {g_D}_j\left(-\underbrace{\left(\varphi_j,\pad{\varphi_i}{n}\right)_{L^2(\Gamma_{D,h})}}_{\rm D_{\Gamma_D}} \!\!\!\!\!\!\!+ \underbrace{\lambda(\varphi_j,\varphi_i)_{L^2(\Gamma_{D,h})}}_{\rm P_{\Gamma_D}} \right)
\end{align}
This result is the discrete system $A^{\rm FEM} u_h = F^{\rm FEM}$,
where 
\begin{align}
\label{eq:LS_FEM}
A^{\rm FEM} = {\rm S - S_T} +\lambda {\rm P_{\Gamma_D}},\\ \nonumber
F^{\rm FEM} = {\rm M} f + (\lambda{\rm P_{\Gamma_D}- D_{\Gamma_D}}) g_D + {\rm N_{\Gamma_N}}\, g_N.
\end{align}

\section{Further remarks}
\label{sec:remarks}
In this section, we would like to comment on other aspects that differentiate the $\mathghost$-FEM and the Coco-Russo scheme. In particular, we will focus our attention on the theoretical results regarding the convergence rate of the here proposed schemes.
We will also discuss the \textit{small-cut} issue arising in $\mathghost$-FEM and the availability of fast-solvers for the proposed schemes.
\subsection{A priori convergence estimates}
The convergence analysis of the Coco-Russo scheme is still an open problem.
In \cite{cocoStissi}, using the Toeplitz operator theory and generalized locally Toeplitz matrix-sequences, it has been proven that in the one-dimensional case the Coco-Russo scheme converges with second-order accuracy with respect to the ${L}^\infty$ norm.
In fact, the authors have shown that in one dimension the $\lVert \left(A^{\rm FDM}\right)^{-1}\lVert_\infty$ remains constant as the mesh size $h$ goes to zero, which in combination with standard finite difference consistency estimates guarantees second-order accuracy of the scheme.
Furthermore, in the same work, the authors have partially extended the results to the two-dimensional case, showing that the scheme converges with second-order accuracy with respect to the ${L}^\infty$ norm, if the boundaries of the domain are parallel to the axes.

Taking advantage of the variational formulation underlying the $\mathghost$-FEM, much more general convergence results can be obtained for this scheme. In \cite{astuto2024nodal}, the authors have proven, for general domains and boundary conditions, the convergence of the $\mathghost$-FEM with second-order accuracy with respect to the ${L}^2$ norm.
It is worth mentioning that the convergence analysis of the $\mathghost$-FEM is almost identical to the one of other classical unfitted FEMs, see for example \cite{Burman2015,Lehrenfeld2016}.
Furthermore, for the $\mathghost$-FEM to be convergent with optimal rate, the penalization term $\lambda$ is assumed to be sufficiently large.
In particular, we fixed $\lambda$ equal to $h^{-\alpha}$, with $\alpha\in [\frac{3}{2},2]$, because a large penalization term might lead to a ill-conditioned linear system.
As we will discuss in Section \ref{sec:results}, keeping a fixed $\lambda$ proportional to the mesh sizes requires us to address the so-called \textit{small-cut} problem.
\subsection{ \emph{Small-cut} cell problem for {$\protect \mathghost$-FEM}}\label{sect:cond}
A well-known issue affecting all unfitted discretization is the so-called \textit{small-cut} problem, which arises when the intersection between the domain $\Omega$, described via a level-set function, and the grid cells is very small.
When resorting to a Nitsche formulation, as we do in the $\mathghost$-FEM, the \textit{small-cut} problem arises. This occurs because, 
in order to preserve the coercivity of the bilinear form, the penalization term $\lambda$ must be of the same order as the inverse of the area of the small cells obtained as a consequence of such \textit{small-cut}.

Different approaches have been proposed to deal with such an issue. Among them, the \textit{ghost penalty} technique is the defacto preferred one in the CutFEM. The idea behind the \textit{ghost penalty} is to add a penalization term to the bilinear form, which weakly enforces higher continuity across element interfaces coupling basis functions with small support to larger neighboring elements, \cite{Burman2015}.
Another approach is based on \textit{cell agglomeration}, which consists of merging small cells with their neighbors in order to obtain larger cells, \cite{aggregation}.

Lastly, we point out that the simplest approach would be to exclude the basis functions associated with small cells from the system, \cite{snapping}. In \cite{astuto2024nodal}, 
this approach is viewed as a perturbation of the original domain and it is proven that the $\mathghost$-FEM is convergent with optimal rate even in the presence of small cells. In practise, 
we evaluate the level set function $\phi$ at the vertices of each cell: if the value is {smaller than a threshold} equal to a power of the length of the cell, i.e.
{if $0<-\phi<h^\alpha$, where $\alpha \in [\frac{3}{2},2]$ is the snapping exponent,} we disregard the respective cell.

\begin{figure}
\centering
\begin{minipage}{0.49\textwidth}
\centering
\begin{overpic}[abs,width=0.75\textwidth,unit=1mm,scale=.25]{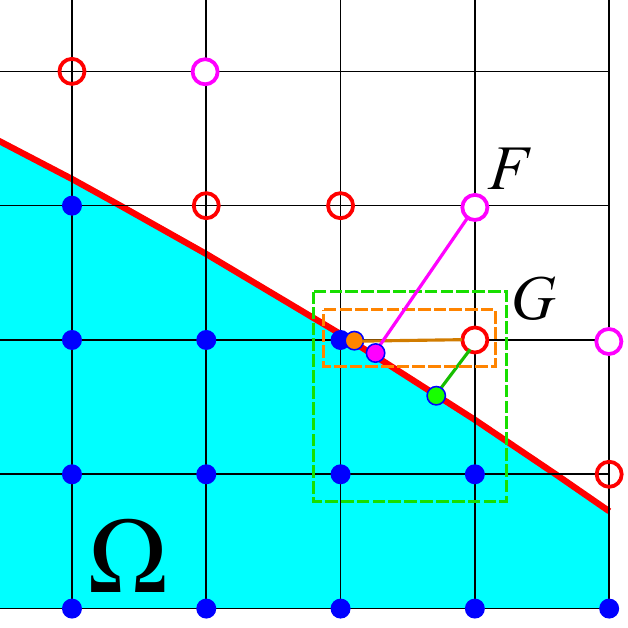}
\put(-2,45){(a)}
\end{overpic}
\end{minipage}
\begin{minipage}{0.49\textwidth}
\centering
\begin{overpic}[abs,width=0.75\textwidth,unit=1mm,scale=.25]{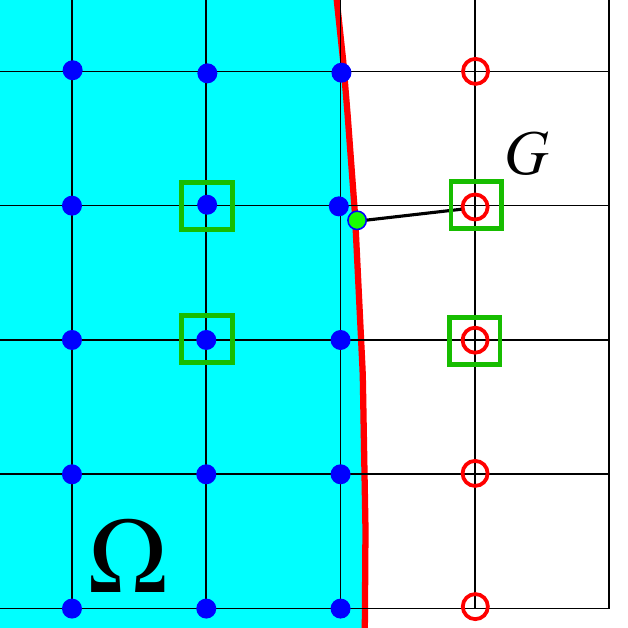}
\put(-2,45){(b)}
\end{overpic}
\end{minipage}
\caption{\textit{(a): Ghost value extrapolation would be ill-conditioned for the external grid point $F$. However, $F$ is not a ghost point of the finite-difference method. The extrapolation of the ghost value $G$ would be ill-conditioned if the boundary condition is enforced to the orange boundary point and a linear interpolation stencil (orange dashed rectangle) is used, as in~\cite{Gibou2002}. The orthogonal projection and the bilinear interpolation stencil (green dashed rectangle) adopted in the Coco-Russo method improves the conditioning. (b): The ill-conditioned extrapolation of the ghost value $G$ is overcome by enlarging the interpolation stencil (green squares).}}
\label{fig:small_cell_FDM}
\end{figure}

{
\subsection{Ill-conditioned ghost value extrapolation for \\ finite-difference methods}
The equivalent of the small cell problem in finite-difference methods is the ill-conditioned behaviour that arises when the diagonal entries of the linear equations \eqref{QHghost} or \eqref{QHghost2} are small compared to the off-diagonal entries.

From \eqref{coeffsLSstencil} to \eqref{eq:normal_deriv}, the diagonal entry is $l_0(\vartheta_x) l_0(\vartheta_y)$ for Dirichlet boundary conditions and $- \left( |n_x| l'_0 (\vartheta_x) l_0(\vartheta_y) + |n_y| l_0 (\vartheta_x) l_0'(\vartheta_y) \right)$ for Neumann boundary conditions. From \eqref{coeffs_l1} and \eqref{coeffs_l2}, considering that $0 \leq \vartheta_x, \vartheta_y < 1$, The issue can be observed only for Dirichlet boundary conditions and when at least one of $\vartheta_x$ or $\vartheta_y$ is close to one. From Fig.~\ref{fig:small_cell_FDM} (a), we observe that these effects are less prevalent than small cell problems in $\mathghost$-FEM. In fact, the effect would appear for purple external points, such as $F$.
However, purple points are external grid points in $\mathghost$-FEM but they are not ghost points in the finite-difference method, and thus the effect is not observed. Now consider the ghost point \(G\), which is a ghost point in the finite-difference method. Some ghost point approaches, such as~\cite{Gibou2002}, enforce the boundary condition on the horizontal or vertical boundary projection. This approach can be problematic in some cases: for example, the horizontal boundary projection of $G$ in the panel (a) of Fig.~\ref{fig:small_cell_FDM} (orange point) leads to $\vartheta_x \approx 1$. The Coco-Russo method uses the orthogonal projection instead and this drastically reduces the likelihood of ill-conditioned effects.

Finally, there are scenarios where the issue does affect the Coco-Russo method, as shown in the (b) panel of Fig.~\ref{fig:small_cell_FDM}. In such cases, we overcome the problem by enlarging the interpolation stencil. For example, in the case of bilinear interpolation (\(p=1\)), we use the stencil depicted with green squares. This adjustment effectively halves the value of $\vartheta_x$ in the interpolation formula. This modification can be applied when $\left|1-\vartheta_x\right| < \varepsilon$, where $\varepsilon>0$ is a suitable tolerance, such as $\varepsilon = h$. The same argument applies to $\vartheta_y$. We have observed numerically that this modification is necessary only in very rare cases.
}

\subsection{Fast solvers}
{
Linear systems resulting from the finite-difference discretization of elliptic equations are typically symmetric and positive definite, making them amenable to standard fast solvers such as conjugate gradient methods. However, the ghost-point approach in the Coco-Russo method compromises both symmetry and positive definiteness. In \cite{Coco2013}, the authors demonstrated that the Gauss-Seidel method fails to converge for the Coco-Russo method. To address this, they implemented a relaxation strategy for the boundary conditions, wherein the equations associated with ghost points are relaxed using a suitable parameter to ensure convergence.

To further enhance convergence, a tailored geometric multigrid solver is employed, specifically designed to handle curved boundaries. The transfer operators (restriction and interpolation) are designed to separately process internal and ghost equations, preventing their contributions from mixing during grid transfers between different levels of resolution. For more details, refer to \cite{Coco2013}.

In multigrid methods for complex-shaped domains, boundary effects can significantly impact overall efficiency unless properly addressed. To mitigate this, it is common practice to add extra relaxations on boundary conditions, maintaining overall performance~\cite{Brandt1984}. This strategy, used in~\cite{Coco2013,Coco2018,COCO2020109623}, involves minimal additional computation compared to internal relaxations, becoming negligible as the spatial step $h \rightarrow 0$ \cite{Trottemberg:MG}. While effective for uniformly resolved grids, real-world problems often require finer meshes near boundaries to capture curvature and maintain uniform numerical error. In such cases, the cost of extra boundary relaxations can dominate, leading to a sub-optimal solvers.

In~\cite{coco2023ghost} the authors developed a Boundary Local Fourier Analysis (BLFA) theory for the Coco-Russo methods to design an efficient relaxation scheme that smooths the residual of boundary conditions along the tangential direction without compromising the smoothing performance of the internal equations. This smoother is integrated into a multigrid framework, ensuring that the convergence factor remains unaffected by boundary effects.

It is well-known in the literature that preconditioning a positive-definite linear system is an easier task than dealing with an indefinite system. In Table \ref{tab:precond} we explore different preconditioners for the symmetric positive-definite linear system originating from the $\mathghost$-FEM discretization. In particular, we explore the Jacobi and Successive Over Relaxation preconditioners implemented in PETSc \cite{petsc} together with the HYPRE \cite{hypre} implementation of Algebraic MultiGrid (AMG) wrapped in PETSc.
We notice from Table \ref{tab:precond}, that even if algebraic multigrid preconditioners outperform all other choices of preconditioners, still the number of conjugate gradient iterations depends on the number of degrees of freedom. For this, we plan to further investigate the development fast solvers for the $\mathghost$-FEM. 
\begin{table}
    \caption{We compare different preconditioners applied to a $\protect\mathghost$-FEM discretization of the Poisson equation formulated on a circle, with an increasing number of degrees of freedom. We display the number of conjugate gradient iterations required to reach the residual denoted between parenthesis when using no precondition, Jacobi preconditioner, Successive Over Relaxation preconditioner and HYPRE implementation of Algebraic MultiGrid preconditioners.}
    \label{tab:precond}
    \begin{tabular}{ccccc} \toprule
    {Degrees of freedom} & {Identity} & {Jacobi} & {SOR} & {AMG (Hypre)} \\ \midrule
    263169  & 991 (1.56e-13) & 430 (2.41e-12)\;\; & 267 (1.32-12)\;\;  & 76 (4.67e-13)  \\
    
    1050625  & 10000 (2.64e-05)\;\; & 839 (3.16e-12)\;\; & 482 (2.04e-12)\;\; & 102 (1.56e-12) \\

    4198401  & 10000 (4.96e-04)\;\; & 1623 (4.27e-12)\;\; & 904 (1.97e-12)\;\; & 151 (7.77e-13) \\
\end{tabular}
\end{table}
}

\section{Numerical experiments}\label{sect:results}
\label{sec:results}
In this section, we compare the results obtained with the two different numerical schemes described in the previous sections. We consider domains $\Omega \subset R = [-1,1]^2$ whose shape is implicitly known through a level-set function. We calculate the relative error between the numerical and exact solutions, and between the gradient of the numerical solutions and of the exact one, with the following formulas 
\begin{equation}
\label{eq:error}
    {\rm error}  = \frac{||{\tt f}_{h} - {\tt f}_{\rm exa}||_{L^\beta(\Omega_h)}}{||{\tt f}_{\rm exa}||_{L^\beta(\Omega_h)}}, 
\end{equation}
where $\beta = 1,2,\infty$, ${\tt f} = u$ for the solution and ${\tt f} = \nabla u$ for its gradient.
\begin{figure}[H]
\hspace{-0.5cm}
	\centering
	\begin{minipage}
		{.24\textwidth}	\centering
\begin{overpic}[abs,width=1.2\textwidth,unit=1mm,scale=.25]{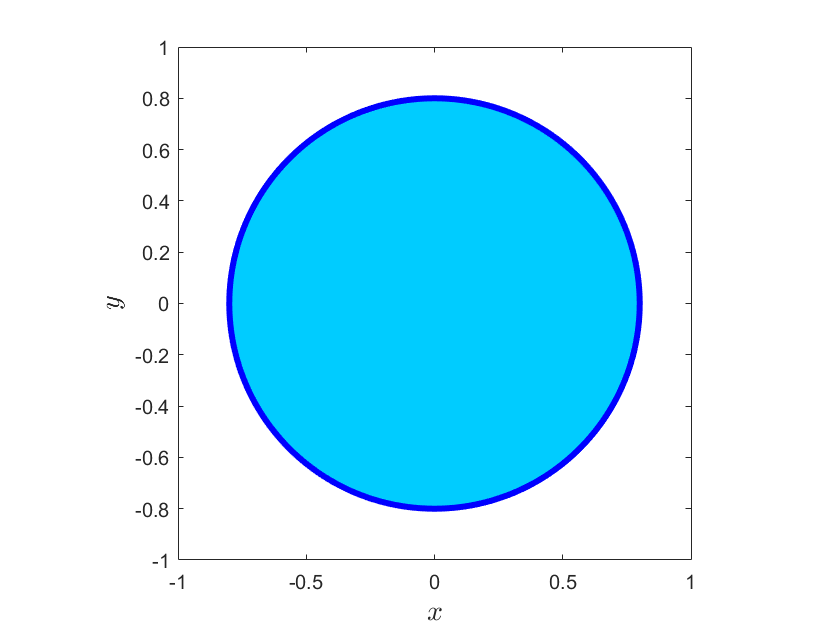}
\put(2.,26){(a)}
\end{overpic}
	\end{minipage}
	\begin{minipage}
		{.24\textwidth}	\centering
\begin{overpic}[abs,width=1.2\textwidth,unit=1mm,scale=.25]{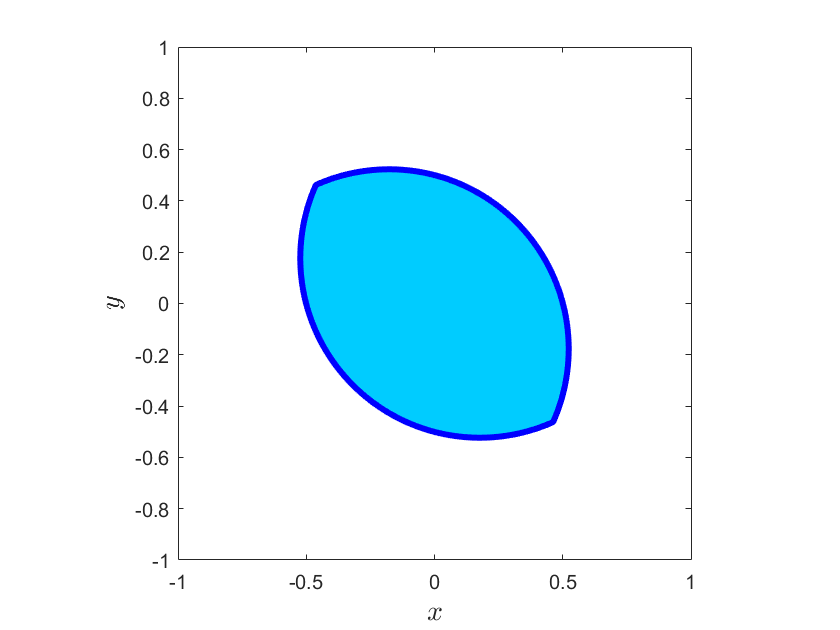}
\put(2.,26){(b)}
\end{overpic}
	\end{minipage}
	\begin{minipage}
		{.24\textwidth}	\centering
\begin{overpic}[abs,width=1.2\textwidth,unit=1mm,scale=.25]{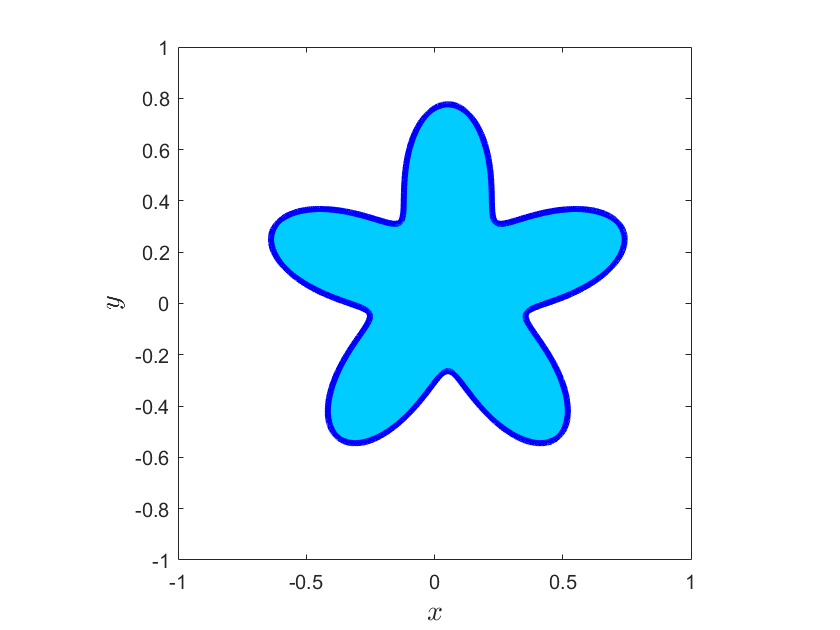}
\put(2.,26){(c)}
\end{overpic}
	\end{minipage}
	\begin{minipage}{.24\textwidth}	
  \centering
\begin{overpic}[abs,width=1.2\textwidth,unit=1mm,scale=.25]{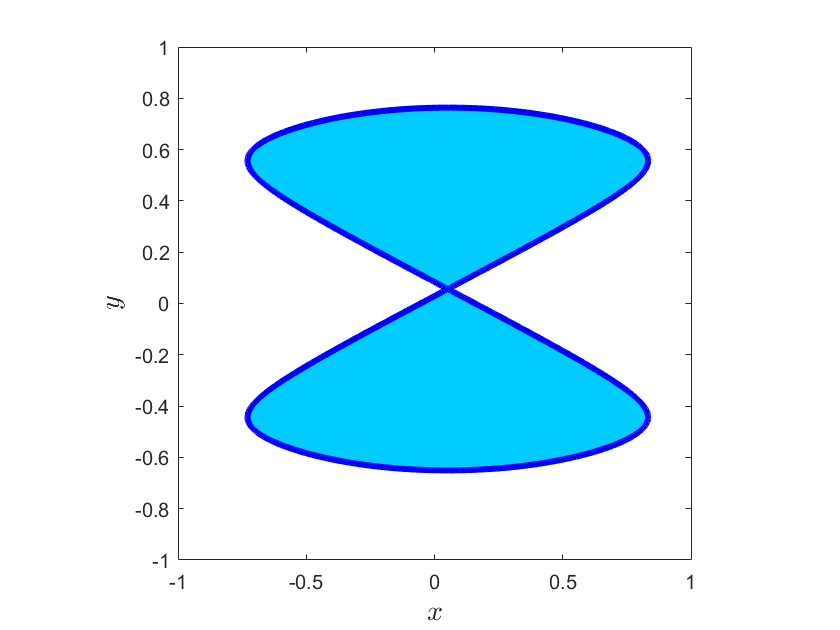}
\put(2.,26){(d)}
\end{overpic}
	\end{minipage}
	\caption{\textit{Shape of the domains considered in our tests. We have the circular (a), the rotated leaf- (b), the flower- (c) and the hourglass-shaped domain (d).}}
	\label{fig:shape_domains}
\end{figure}
We also compare the behaviour of the conditioning number ${\rm cond}(A^{\rm FDM})$ and ${\rm cond}(A^{\rm FEM})$ of the linear systems in \eqref{eq:LS_FDM},(\ref{QHghost}--\ref{QHghost2}) and \eqref{eq:LS_FEM}, respectively. We conduct a deeper investigation into the linear system described in \eqref{eq:LS_FEM}, presenting the results obtained when varying the snapping exponent 
$\alpha$ in the penalization term $\lambda (= h^{-\alpha})$.

We show the numerical results with Dirichlet and mixed boundary conditions. In our tests, we use different expressions for level-set functions to define various domains (see Fig.~\ref{fig:shape_domains}), while the exact solution is 
\[ u_{\rm exa} = \sin(x)\sin(y).
\]
\subsubsection*{Circular domain} 
\begin{figure}[H]\hspace{-1.3cm}
	\centering
	\begin{minipage}
		{.32\textwidth}	\centering
\begin{overpic}[abs,width=1.3\textwidth,unit=1mm,scale=.25]{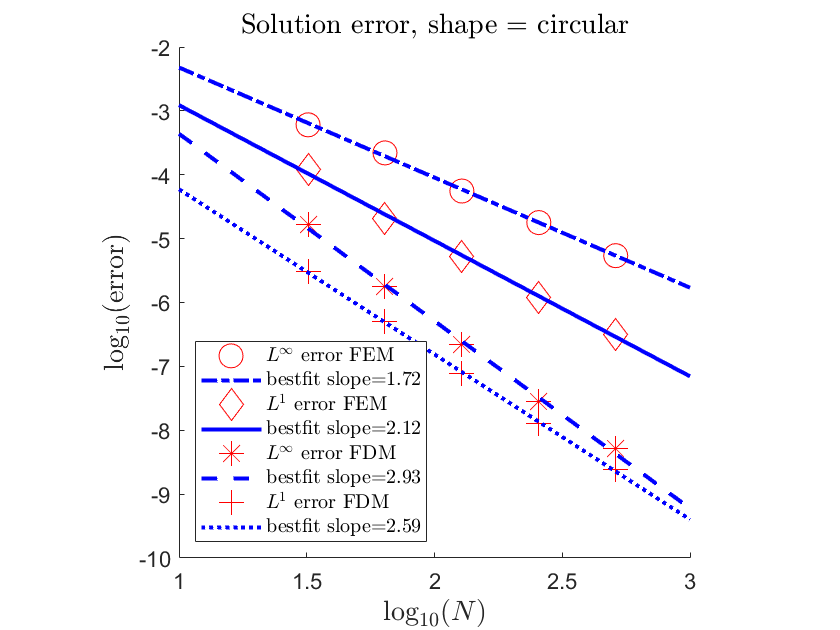}
\put(4,37){(a)}
\end{overpic}
	\end{minipage}
	\begin{minipage}
		{.32\textwidth}	\centering
\begin{overpic}[abs,width=1.3\textwidth,unit=1mm,scale=.25]{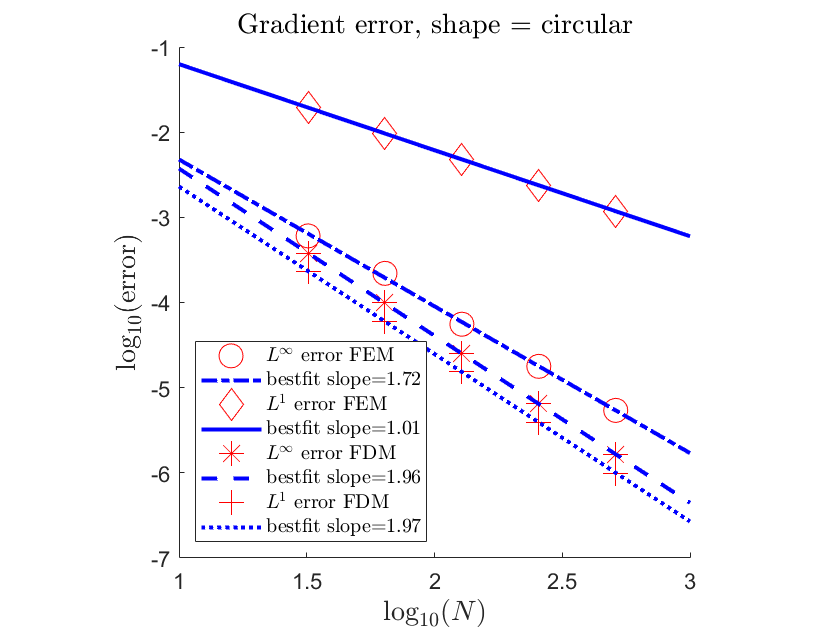}
\put(4,37){(b)}
\end{overpic}
	\end{minipage}
	\begin{minipage}
		{.32\textwidth}	\centering
\begin{overpic}[abs,width=1.3\textwidth,unit=1mm,scale=.25]{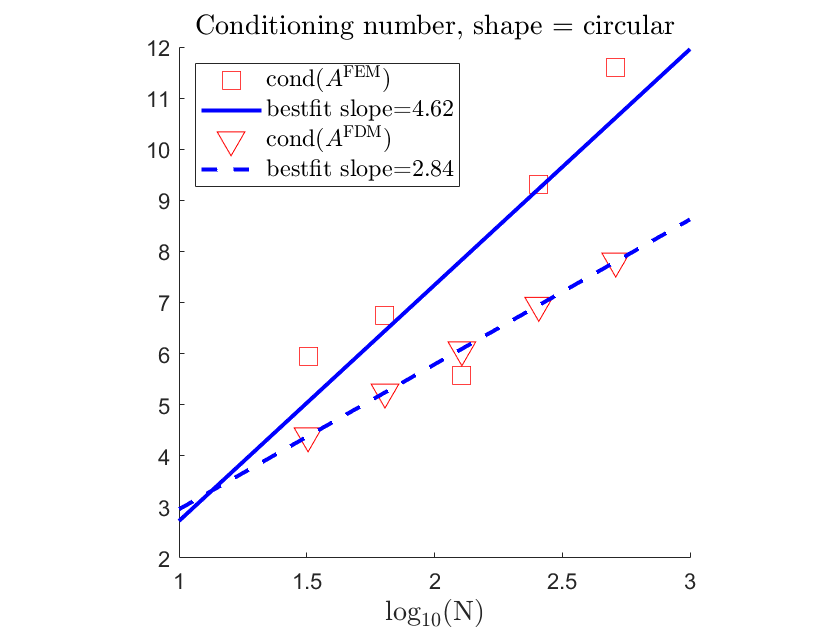}
\put(4,37){(c)}
\end{overpic}
	\end{minipage}
	\caption{\textit{Comparison of the error behavior and of the conditioning number between the two different numerical methods, for the circular domain and Dirichlet boundary conditions: relative error of the numerical solutions (a), of the gradient (b) and conditioning number of the linear systems in (c); snapping exponent $\alpha = 2$.}}
	\label{fig:results_circle}
\end{figure}

Let us start with a circular domain  centered at the origin $(0,0) \in \Omega$. The level-set function is $\phi=r-\sqrt{x^2+y^2}$, where $r = 0.8$ is the radius of the circle. 
We consider Dirichlet boundary conditions (i.e.  $\Gamma_D = \Gamma$) in Fig.~\ref{fig:results_circle}, and mixed boundary conditions in Fig.~\ref{fig:results_circle_mixed}, i.e. $\Gamma_D = \Gamma \cap \{x\leq 0\}$ and $\Gamma_N = \Gamma \cap \{x > 0\}$. 

\begin{figure}[H]\hspace{-1.3cm}
	\centering
	\begin{minipage}
		{.32\textwidth}	\centering
\begin{overpic}[abs,width=1.3\textwidth,unit=1mm,scale=.25]{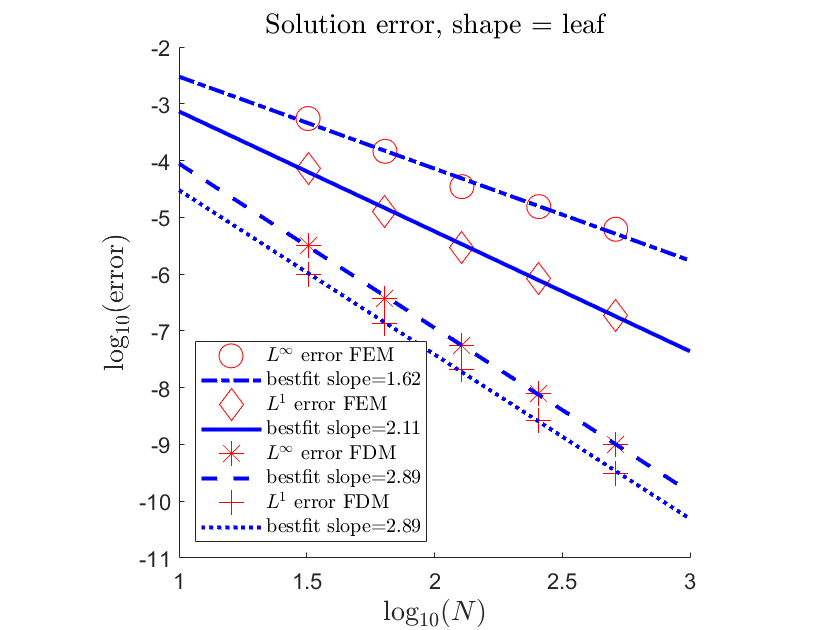}
\put(4,37){(a)}
\end{overpic}
	\end{minipage}
	\begin{minipage}
		{.32\textwidth}	\centering
\begin{overpic}[abs,width=1.3\textwidth,unit=1mm,scale=.25]{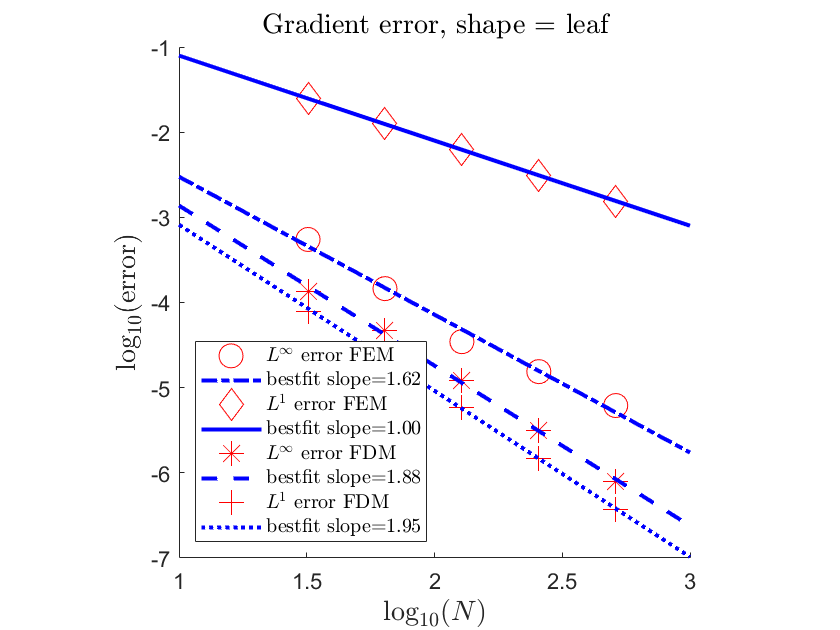}
\put(4,37){(b)}
\end{overpic}
	\end{minipage}
	\begin{minipage}
		{.32\textwidth}	\centering
\begin{overpic}[abs,width=1.3\textwidth,unit=1mm,scale=.25]{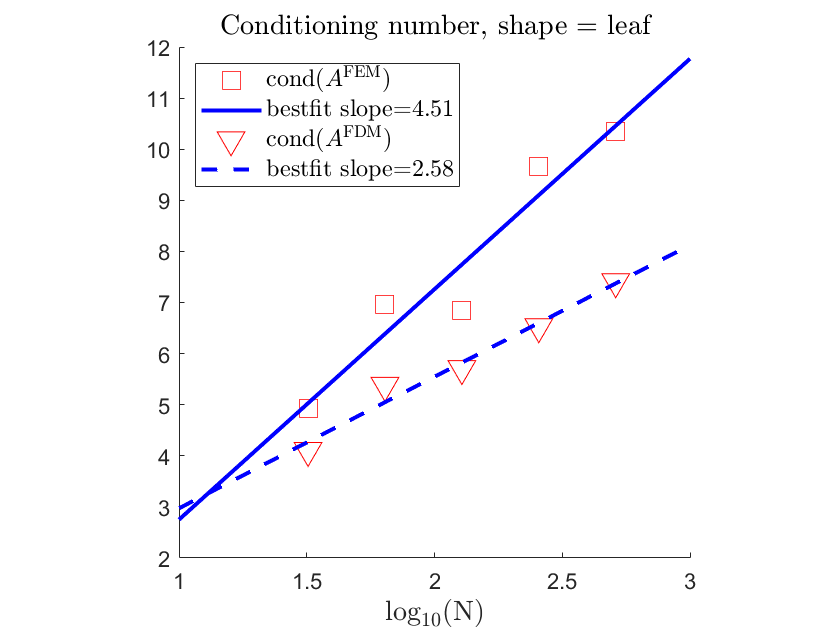}
\put(4,37){(c)}
\end{overpic}
	\end{minipage}
	\caption{\textit{Comparison of the error behavior and of the conditioning number between the two different numerical methods, for the leaf-shaped domain and Dirichlet boundary conditions: relative error of the numerical solutions (a), of the gradient (b) and conditioning number of the linear systems in (c); snapping exponent $\alpha = 2$.}}
	\label{fig:results_leaf}
\end{figure}

\subsubsection*{Leaf-shaped domain} We then consider a leaf-shaped domain, with Dirichlet boundary conditions, i.e. $\Gamma = \Gamma_D$ (see Fig.~\ref{fig:results_leaf}), and  mixed boundary conditions, such that $\Gamma_D = \partial \Omega \cap \{x<0\}$ and $\Gamma_N = \partial \Omega \cap \{x\geq 0\}$ (see Fig.~\ref{fig:results_leaf_mixed}). In this case, the level-set function is 
\begin{align*}
&\widetilde x_1 = -0.25,\, x_1 = \widetilde x_1\cos(\pi/4)  \quad R_1 = \sqrt{(x-x_1)^2+y^2} \\
&\widetilde x_2 = 0.25, \,x_1 = \widetilde x_1\sin(\pi/4) \quad R_2 = \sqrt{(x-x_2)^2+y^2} \\         
&r_0 = 0.7,\quad  \phi_1 = R_1-r_0, \quad \phi_2 = R_2-r_0, \quad \phi = \max\{\phi_1,\phi_2\}.
\end{align*}


\subsubsection*{Flower-shaped domain} Here we show a flower-shaped domain, with Dirichlet boundary conditions, i.e. $\Gamma_D = \Gamma$ in Fig.~\ref{fig:results_flower} and mixed boundary conditions in Fig.~\ref{fig:results_flower_mixed}. In this case, the level-set function is 
\begin{align*}
&X = x-0.03\sqrt 3, \quad  Y = y-0.04\sqrt 2, \quad R = \sqrt{X^2+Y^2} \\
&\phi = \frac{R - 0.52 - (Y^5 +5X^4Y-10X^2Y^3)}{5R^5}.   
\end{align*}

\subsubsection*{Hourglass-shaped domain} 
\begin{figure}[H]\hspace{-1.3cm}
	\centering
	\begin{minipage}
		{.32\textwidth}	\centering
\begin{overpic}[abs,width=1.3\textwidth,unit=1mm,scale=.25]{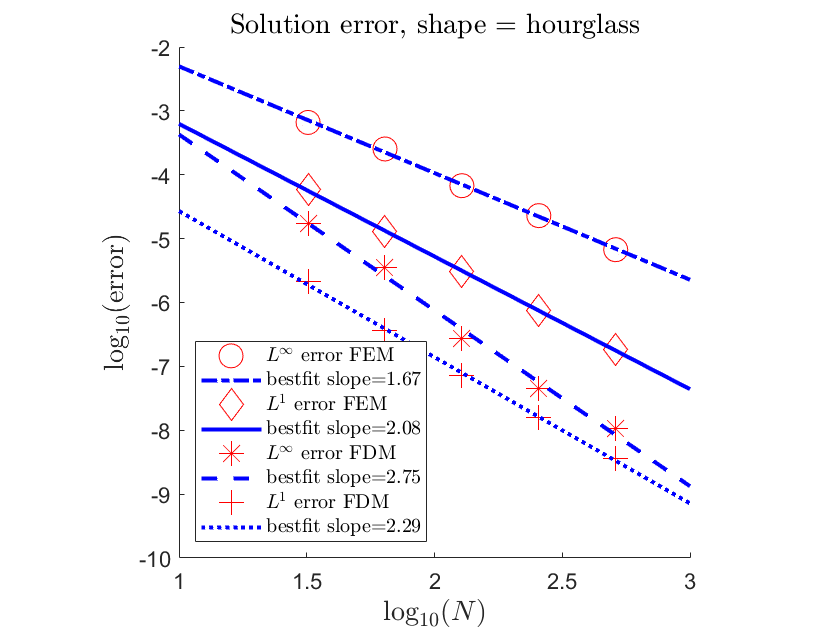}
\put(4,37){(a)}
\end{overpic}
	\end{minipage}
	\begin{minipage}
		{.32\textwidth}	\centering
\begin{overpic}[abs,width=1.3\textwidth,unit=1mm,scale=.25]{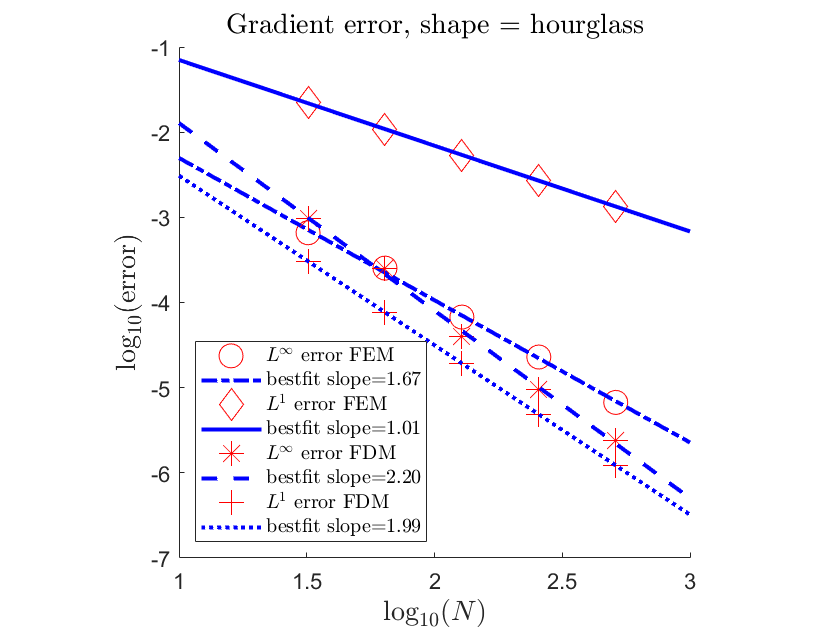}
\put(4,37){(b)}
\end{overpic}
	\end{minipage}
	\begin{minipage}
		{.32\textwidth}	\centering
\begin{overpic}[abs,width=1.3\textwidth,unit=1mm,scale=.25]{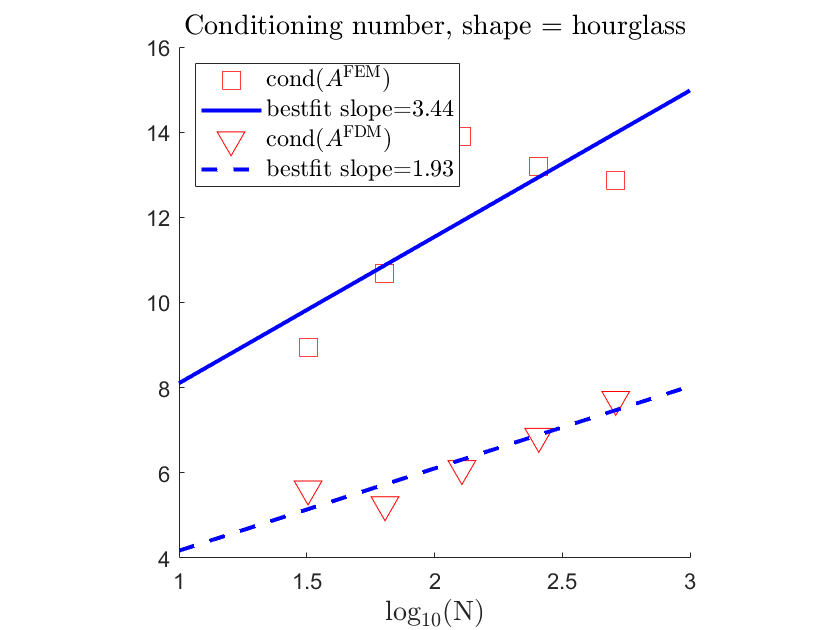}
\put(4,37){(c)}
\end{overpic}
	\end{minipage}
	\caption{\textit{Comparison of the error behavior and of the conditioning number between the two different numerical methods, for the hourglass-shaped domain and Dirichlet boundary conditions: relative error of the numerical solutions (a), of the gradient (b) and conditioning number of the linear systems in (c); snapping exponent $\alpha = 2$.}}
	\label{fig:results_hourglass}
\end{figure}

\begin{figure}[H]\hspace{-1.3cm}
	\centering
	\begin{minipage}
		{.32\textwidth}	\centering
\begin{overpic}[abs,width=1.3\textwidth,unit=1mm,scale=.25]{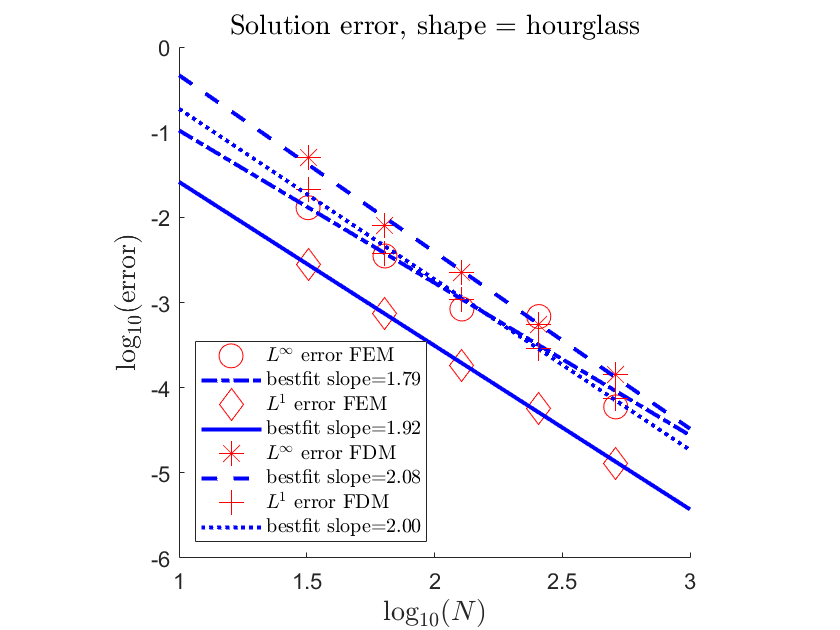}
\put(4,37){(a)}
\end{overpic}
	\end{minipage}
	\begin{minipage}
		{.32\textwidth}	\centering
\begin{overpic}[abs,width=1.3\textwidth,unit=1mm,scale=.25]{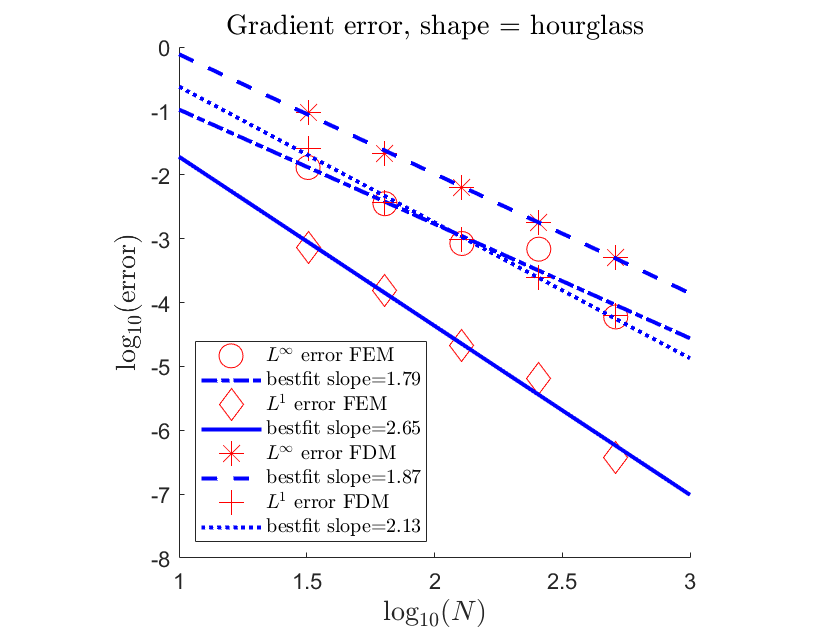}
\put(4,37){(b)}
\end{overpic}
	\end{minipage}
	\begin{minipage}
		{.32\textwidth}	\centering
\begin{overpic}[abs,width=1.3\textwidth,unit=1mm,scale=.25]{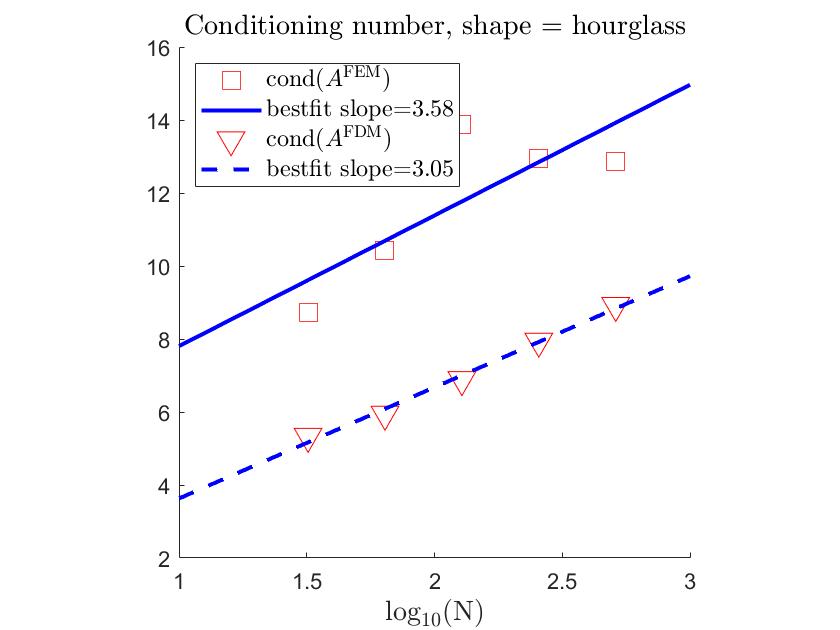}
\put(4,37){(c)}
\end{overpic}
	\end{minipage}
	\caption{\textit{Comparison of the error behavior and of the conditioning number between the two different numerical methods, for the hourglass-shaped domain and mixed boundary conditions: relative error of the numerical solutions (a), of the gradient (b) and conditioning number of the linear systems in (c); snapping exponent $\alpha = 2$.}}
	\label{fig:results_hourglass_mixed}
\end{figure}

Lastly, we have a domain $\Omega$ with a saddle point, whose level-set function is
\begin{align*}
    X = x-0.03\sqrt{3}, \quad Y = y-0.04\sqrt{2}, \quad \phi = 256\,Y^4-16\,X^4-128\,Y^2 + 36\,X^2.
\end{align*}
We choose Dirichlet boundary conditions (i.e.  $\Gamma_D = \Gamma$) in Fig.~\ref{fig:results_hourglass}, and mixed boundary conditions in Fig.~\ref{fig:results_hourglass_mixed}, i.e. $\Gamma_D = \Gamma \cap \{x\leq 0\}$ and $\Gamma_N = \Gamma \cap \{x > 0\}$.

\begin{figure}[H]\hspace{-1.3cm}
	\centering
	\begin{minipage}
		{.32\textwidth}	\centering
\begin{overpic}[abs,width=1.3\textwidth,unit=1mm,scale=.25]{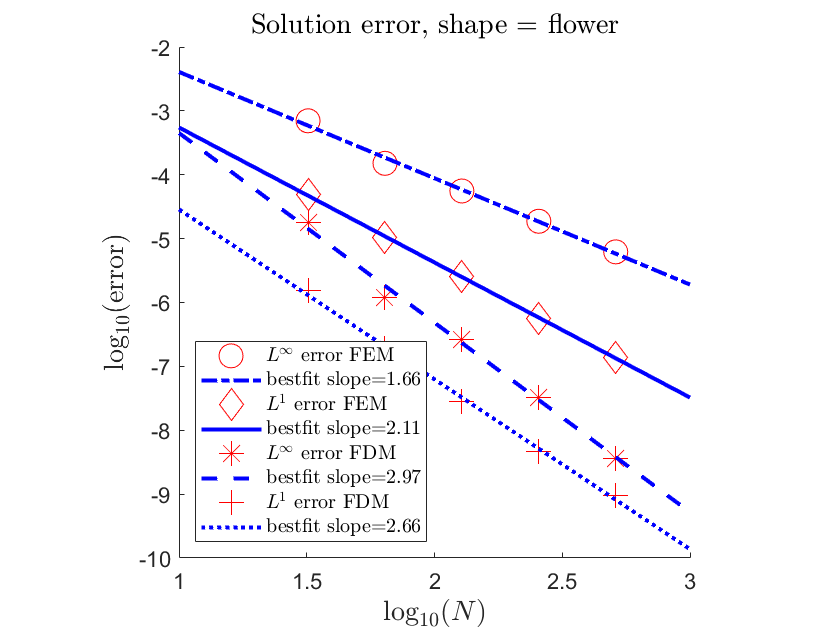}
\put(4,37){(a)}
\end{overpic}
	\end{minipage}
	\begin{minipage}
		{.32\textwidth}	\centering
\begin{overpic}[abs,width=1.3\textwidth,unit=1mm,scale=.25]{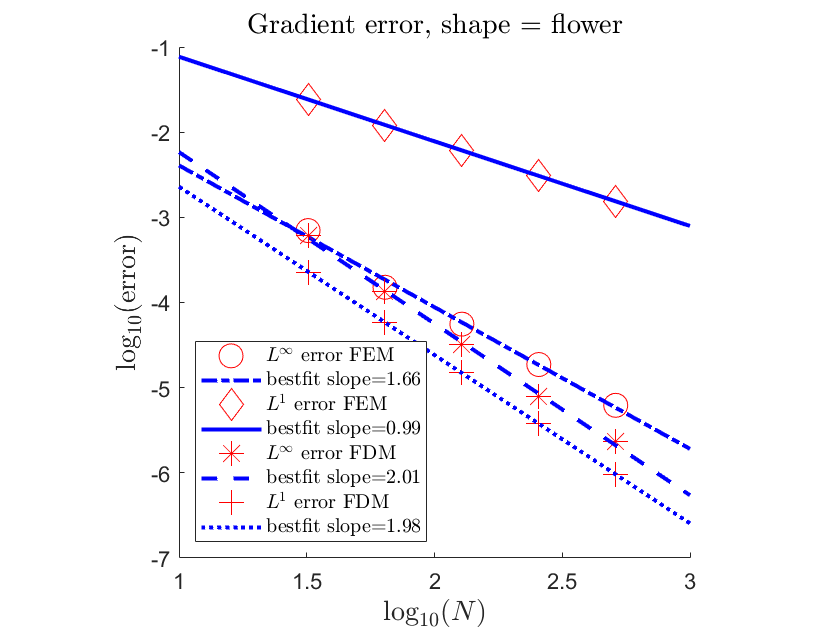}
\put(4,37){(b)}
\end{overpic}
	\end{minipage}
	\begin{minipage}
		{.32\textwidth}	\centering
\begin{overpic}[abs,width=1.3\textwidth,unit=1mm,scale=.25]{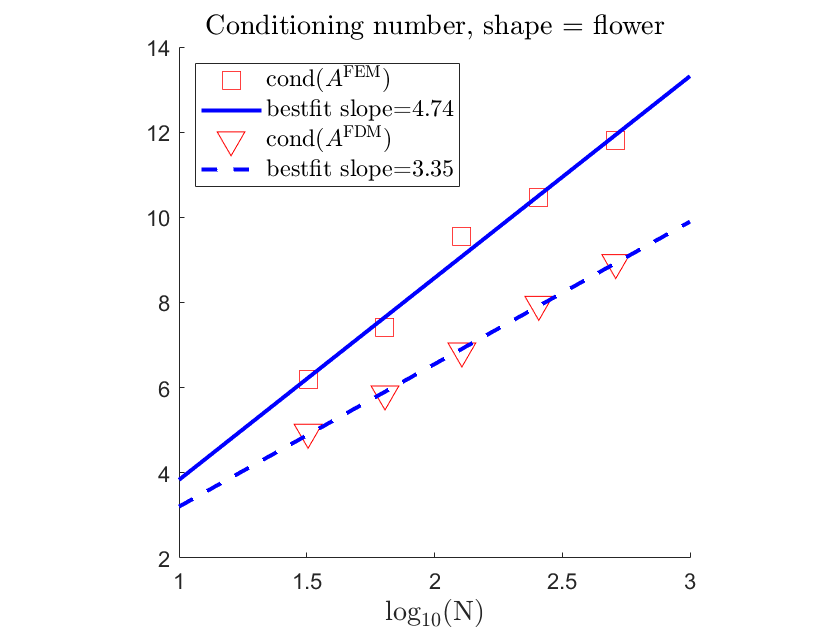}
\put(4,37){(c)}
\end{overpic}
	\end{minipage}
	\caption{\textit{Comparison of the error behavior and of the conditioning number between the two different numerical methods, for the flower-shaped domain and Dirichlet boundary conditions: relative error of the numerical solutions (a), of the gradient (b) and conditioning number of the linear systems in (c); snapping exponent $\alpha = 2$.}}
	\label{fig:results_flower}
\end{figure}

\begin{figure}[H]
\hspace{-1.3cm}
	\centering
	\begin{minipage}
		{.32\textwidth}	\centering
\begin{overpic}[abs,width=1.3\textwidth,unit=1mm,scale=.25]{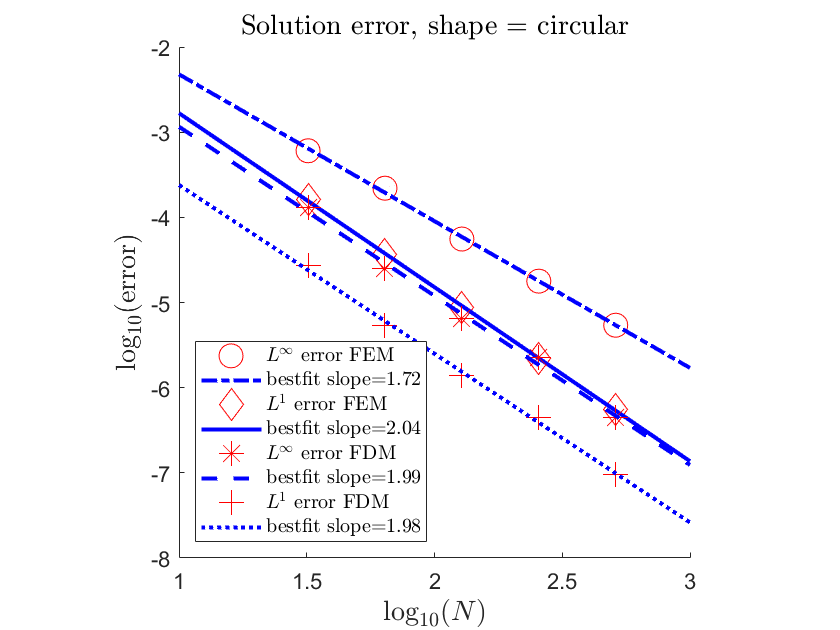}
\put(4.,37){(a)}
\end{overpic}
	\end{minipage}
	\begin{minipage}
		{.32\textwidth}	\centering
\begin{overpic}[abs,width=1.3\textwidth,unit=1mm,scale=.25]{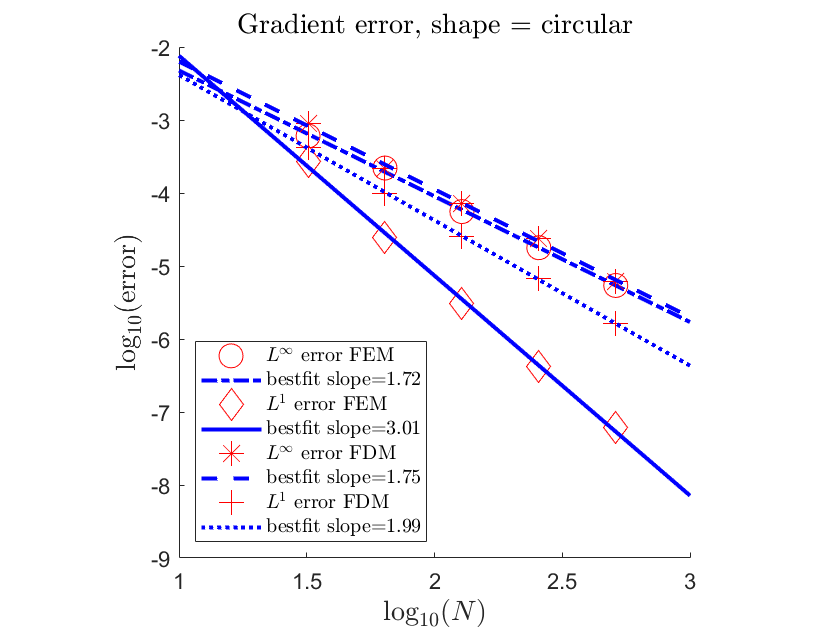}
\put(4.,37){(b)}
\end{overpic}
	\end{minipage}
	\begin{minipage}
		{.32\textwidth}	\centering
\begin{overpic}[abs,width=1.3\textwidth,unit=1mm,scale=.25]{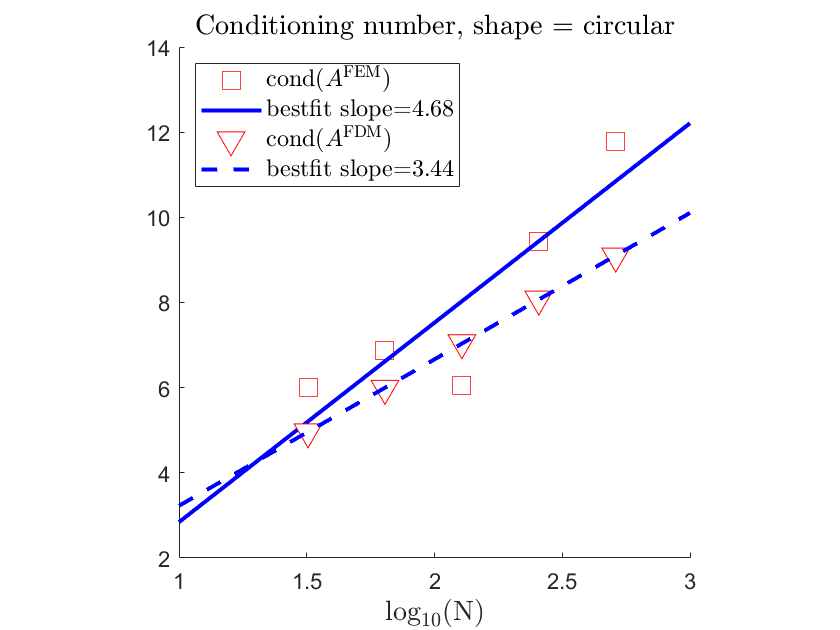}
\put(4,37){(c)}
\end{overpic}
	\end{minipage}
	\caption{\textit{Comparison of the error behavior and of the conditioning number between the two different numerical methods, for the circular domain and mixed boundary conditions: relative error of the numerical solutions (a), of the gradient (b) and conditioning number of the linear systems in (c); snapping exponent $\alpha = 2$.}}
	\label{fig:results_circle_mixed}
\end{figure}

\begin{figure}[H]\hspace{-1.3cm}
	\centering
	\begin{minipage}
		{.32\textwidth}	\centering
\begin{overpic}[abs,width=1.3\textwidth,unit=1mm,scale=.25]{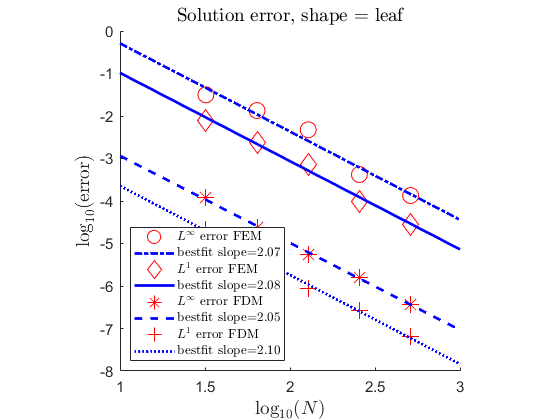}
\put(4,37){(a)}
\end{overpic}
	\end{minipage}
	\begin{minipage}
		{.32\textwidth}	\centering
\begin{overpic}[abs,width=1.3\textwidth,unit=1mm,scale=.25]{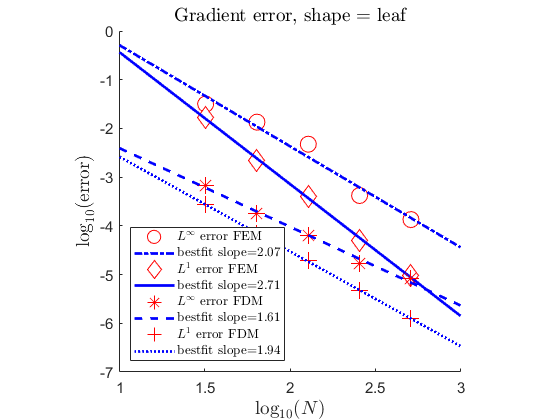}
\put(4,37){(b)}
\end{overpic}
	\end{minipage}
	\begin{minipage}
		{.32\textwidth}	\centering
\begin{overpic}[abs,width=1.3\textwidth,unit=1mm,scale=.25]{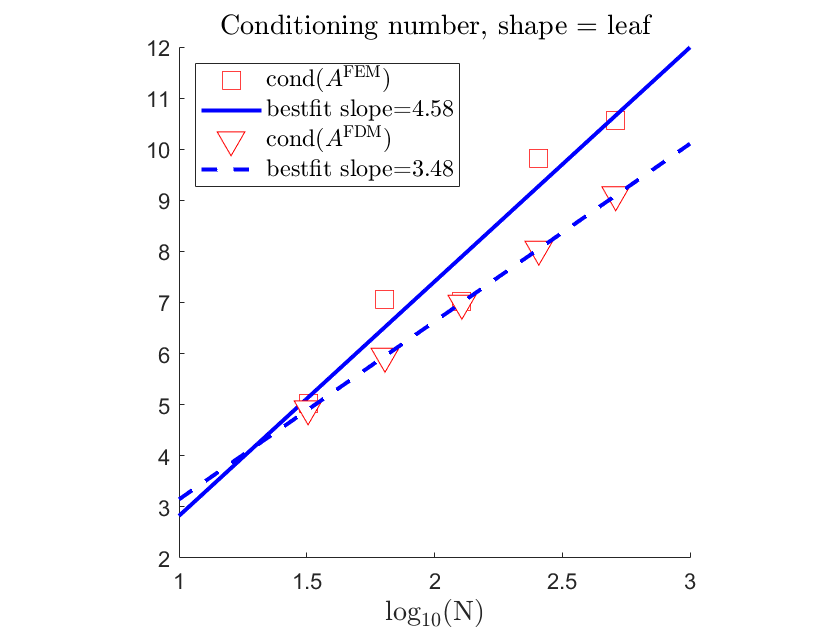}
\put(4,37){(c)}
\end{overpic}
	\end{minipage}
	\caption{\textit{Comparison of the error behavior and of the conditioning number between the two different numerical methods, for the leaf-shaped domain and mixed boundary conditions: relative error of the numerical solutions (a), of the gradient (b) and conditioning number of the linear systems in (c); snapping exponent $\alpha = 2$.}}
	\label{fig:results_leaf_mixed}
\end{figure}

\begin{figure}[htp]\hspace{-1.3cm}
	\centering
	\begin{minipage}
		{.32\textwidth}	\centering
\begin{overpic}[abs,width=1.3\textwidth,unit=1mm,scale=.25]{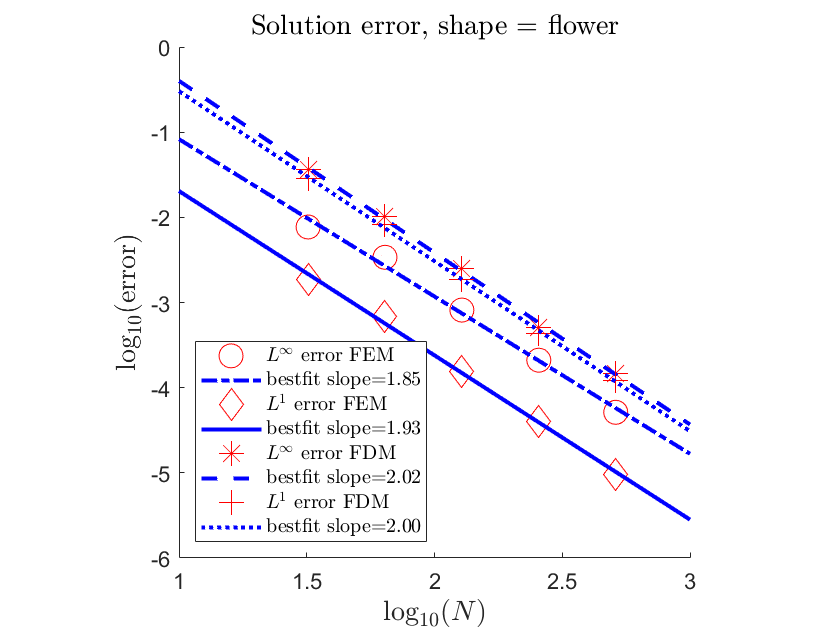}
\put(4,37){(a)}
\end{overpic}
	\end{minipage}
	\begin{minipage}
		{.32\textwidth}	\centering
\begin{overpic}[abs,width=1.3\textwidth,unit=1mm,scale=.25]{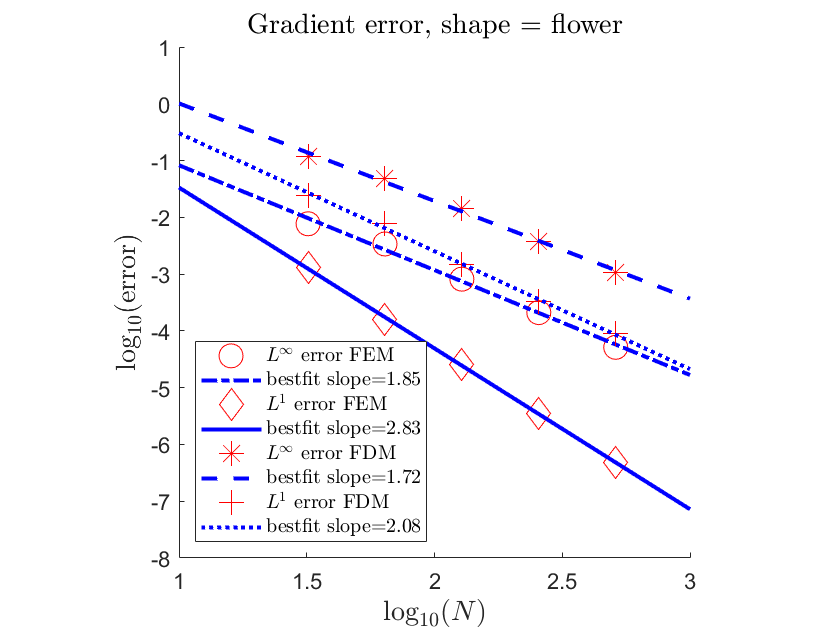}
\put(4,37){(b)}
\end{overpic}
	\end{minipage}
	\begin{minipage}
		{.32\textwidth}	\centering
\begin{overpic}[abs,width=1.3\textwidth,unit=1mm,scale=.25]{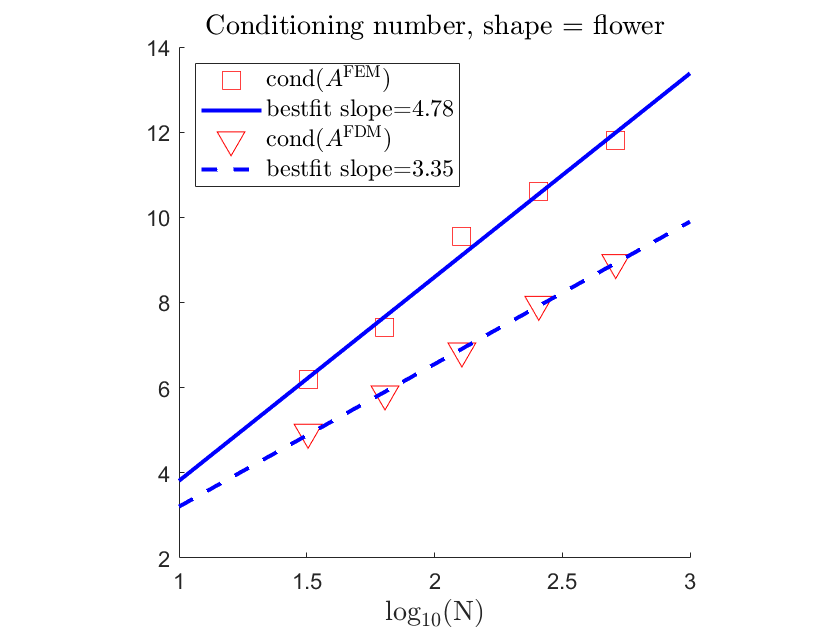}
\put(4,37){(c)}
\end{overpic}
	\end{minipage}
	\caption{\textit{Comparison of the error behavior and of the conditioning number between the two different numerical methods, for the flower-shaped domain and mixed boundary conditions: relative error of the numerical solutions (a), of the gradient (b) and conditioning number of the linear systems in (c); snapping exponent $\alpha = 2$.}}
	\label{fig:results_flower_mixed}
\end{figure}

\begin{figure}[htp]
	\centering
	\begin{minipage}
		{.45\textwidth}	\centering
\begin{overpic}[abs,width=\textwidth,unit=1mm,scale=.25]{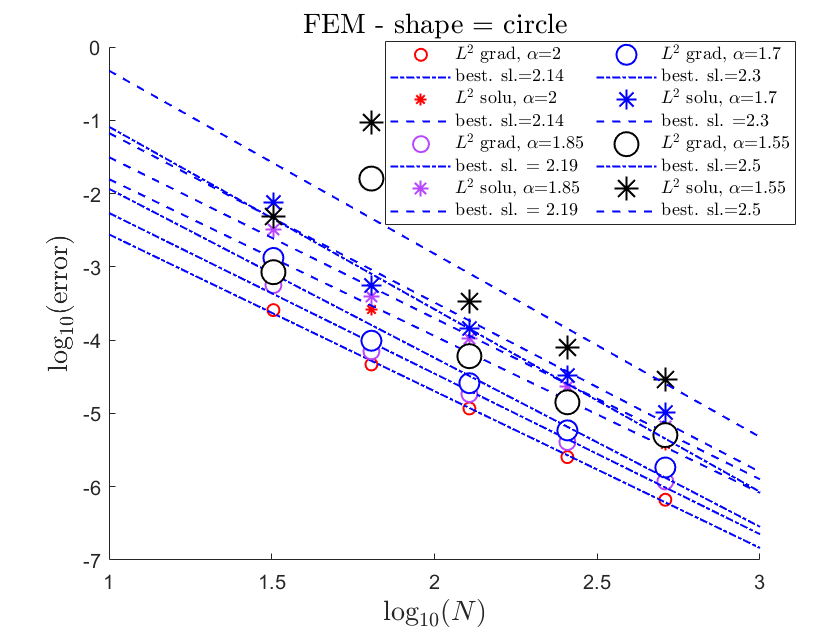}
\put(0.,40){(a)}
\put(16,39){\small $\protect \mathghost$-}
\end{overpic}
	\end{minipage}
	\begin{minipage}
		{.45\textwidth}	\centering
\begin{overpic}[abs,width=\textwidth,unit=1mm,scale=.25]{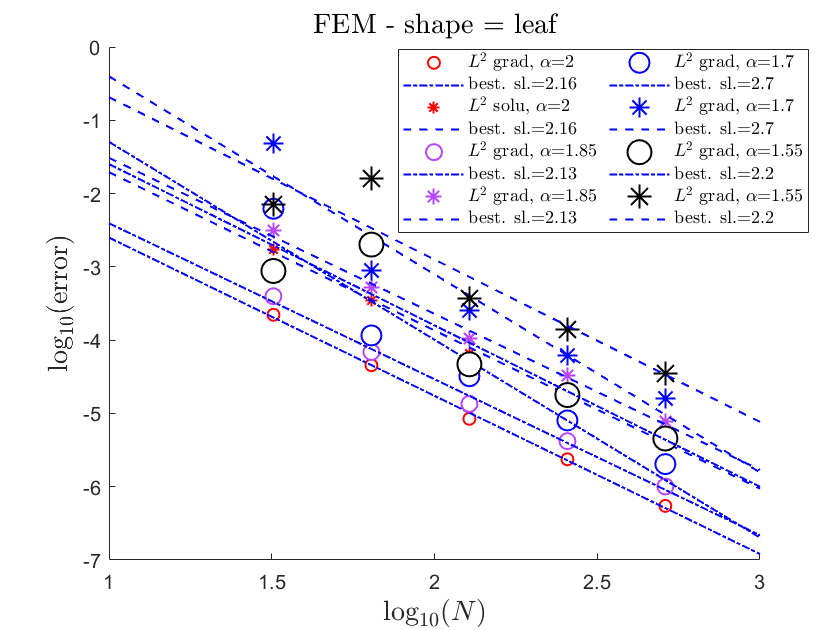}
\put(0.,40){(b)}
\put(16.5,39){\small $\protect \mathghost$-}\end{overpic}
	\end{minipage}
	\begin{minipage}
		{.45\textwidth}	\centering
\begin{overpic}[abs,width=\textwidth,unit=1mm,scale=.25]{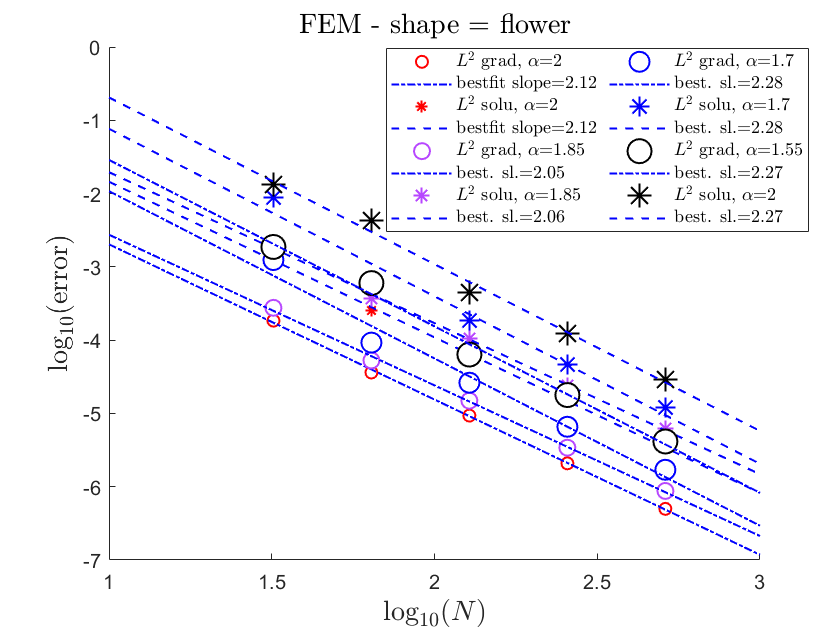}
\put(0.,40){(c)}
\put(16,39){\small $\protect \mathghost$-}\end{overpic}
	\end{minipage}
	\begin{minipage}{.45\textwidth}	
  \centering
\begin{overpic}[abs,width=\textwidth,unit=1mm,scale=.25]{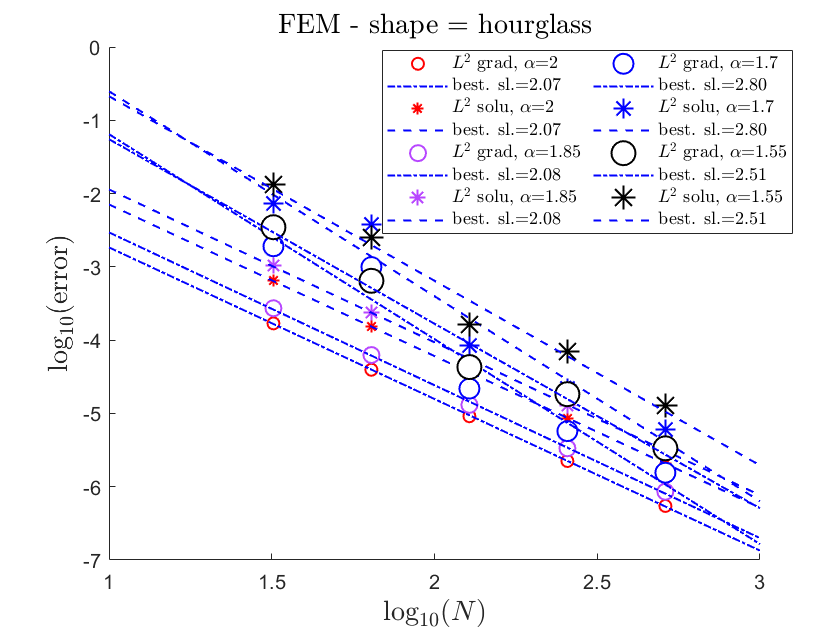}
\put(0.,40){(d)}
\put(14.5,39){\small $\protect \mathghost$-}\end{overpic}
	\end{minipage}
	\caption{\textit{Comparison of the error behaviour changing the snapping exponent $\alpha$ in the snapping threshold and in the penalization term. We show the relative error of the numerical solution and of its gradient. In this test $\alpha = 2, 1.85, 1.7, 1.55$.}}
	\label{fig:error_snapping}
\end{figure}
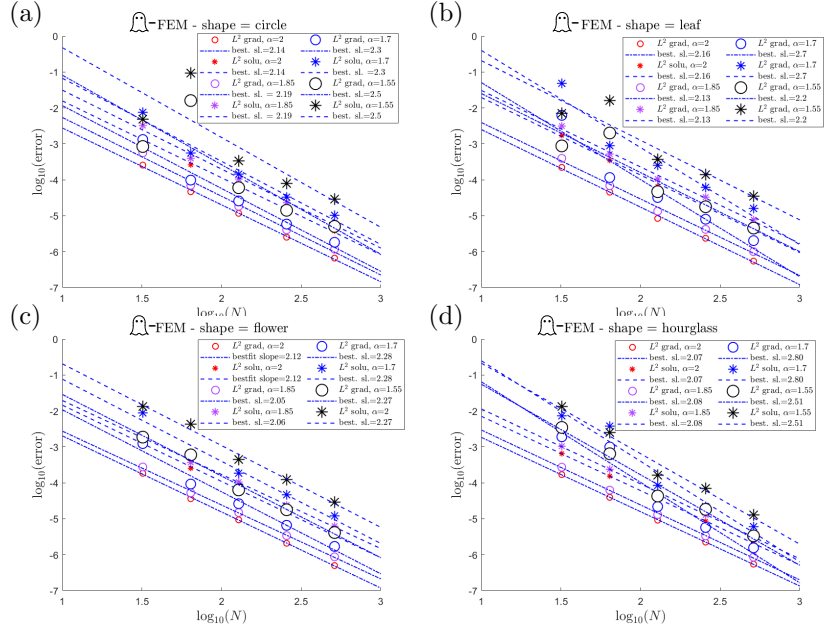

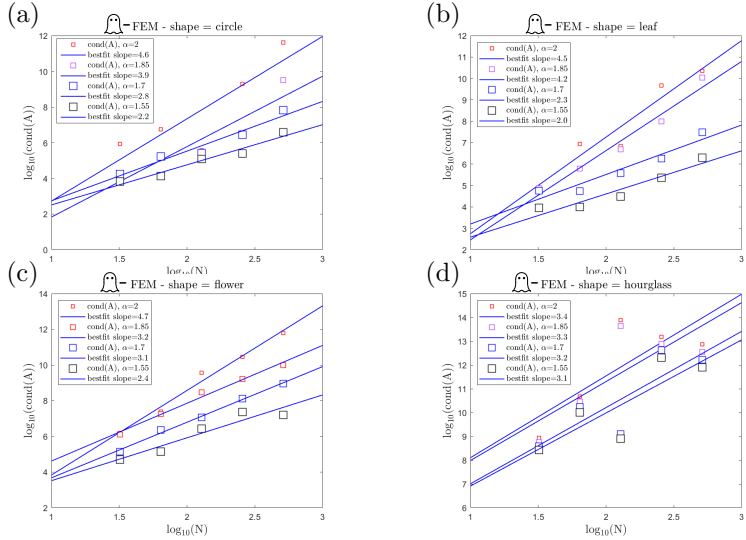
\begin{figure}[htp]
	\centering
	\begin{minipage}
		{.45\textwidth}	\centering
\begin{overpic}[abs,width=0.85\textwidth,unit=1mm,scale=.25]{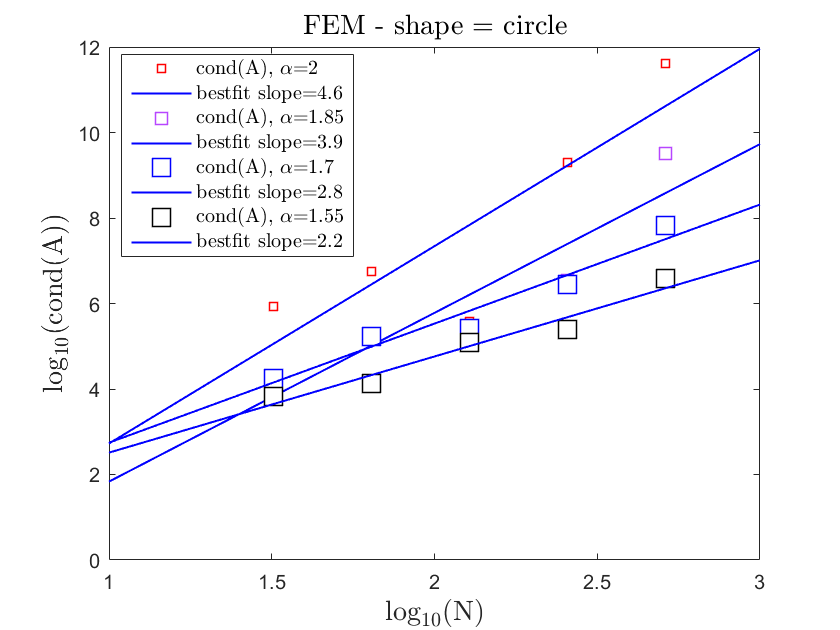}
\put(0.,34){(a)}
\put(13,33){\small $\protect \mathghost$-}\end{overpic}
	\end{minipage}
	\begin{minipage}
		{.45\textwidth}	\centering
\begin{overpic}[abs,width=0.85\textwidth,unit=1mm,scale=.25]{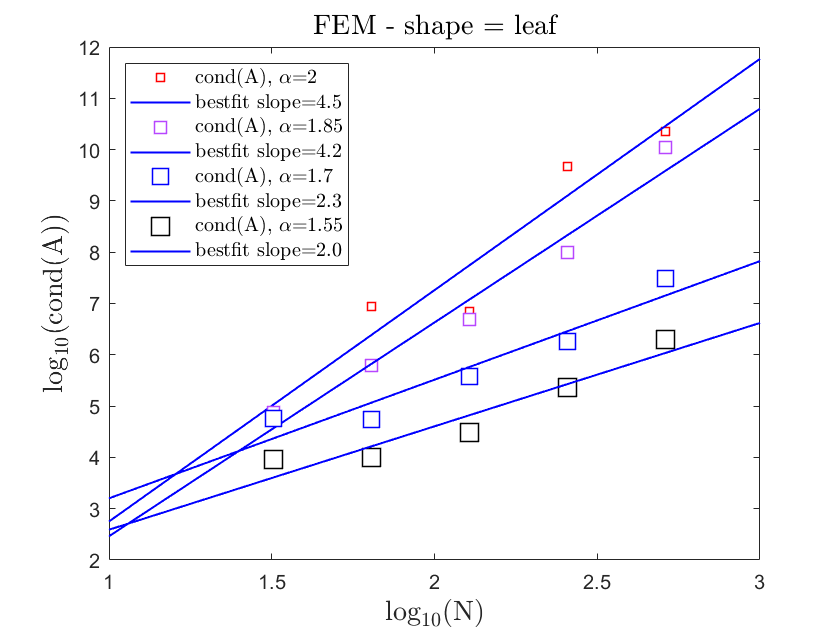}
\put(0.,34){(b)}
\put(13.3,33){\small $\protect \mathghost$-}\end{overpic}
	\end{minipage}
	\begin{minipage}
		{.45\textwidth}	\centering
\begin{overpic}[abs,width=0.85\textwidth,unit=1mm,scale=.25]{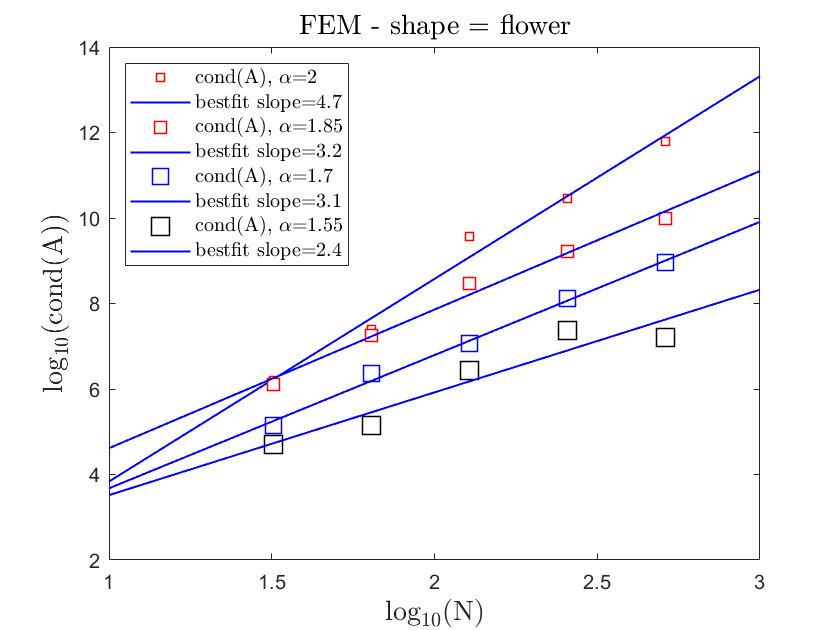}
\put(0.,34){(c)}
\put(12.5,33){\small $\protect \mathghost$-}\end{overpic}
	\end{minipage}
	\begin{minipage}{.45\textwidth}	
  \centering
\begin{overpic}[abs,width=0.85\textwidth,unit=1mm,scale=.25]{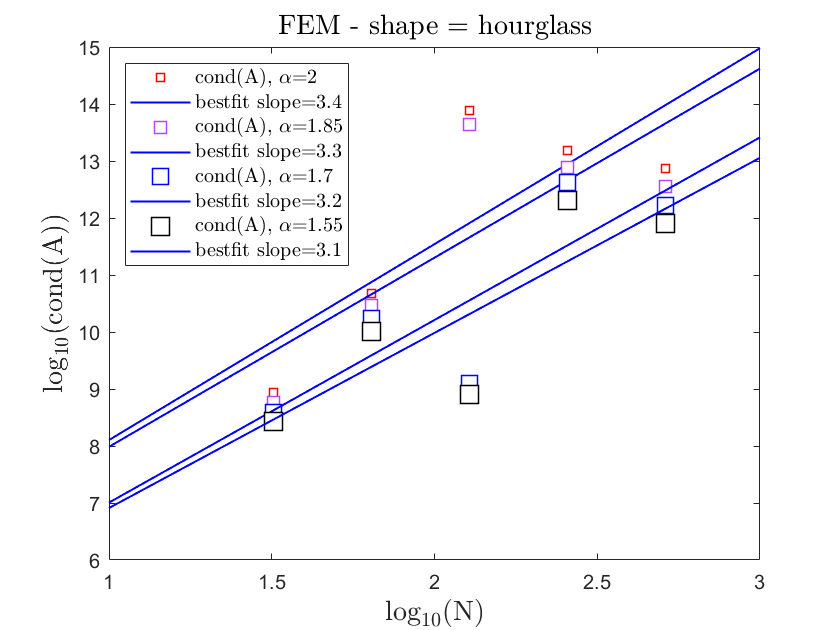}
\put(0.,34){(d)}
\put(11.5,33){\small $\protect \mathghost$-}\end{overpic}
	\end{minipage}
	\caption{\textit{Comparison of the conditioning number behaviour changing the snapping exponent $\alpha$ in the snapping threshold and penalization term. In this test $\alpha = 2, 1.85, 1.7, 1.55$.}}
	\label{fig:cA_snapping}
\end{figure}

We examine various geometries (as illustrated in Fig. \ref{fig:shape_domains}) and types of boundary conditions: in Figs.~\ref{fig:results_circle}--\ref{fig:results_hourglass}, we show the relative error of the solutions (a), of the gradient of the solution (b), and the conditioning number of the linear systems (c), for Dirichlet boundary conditions; analogously, we continue with mixed boundary conditions in Figs.~\ref{fig:results_circle_mixed}--\ref{fig:results_hourglass_mixed}. As previously mentioned, the Coco-Russo scheme employs a 9-point stencil for the interpolation of the solution and of its first space derivatives to ensure a second order accuracy {in the solution and its gradient (a 4-point stencil guarantees second-order accuracy only for Dirichlet boundary conditions and the gradient remains first-order accurate)}. On the contrary, for $\mathghost$-FEM an enlargement of the stencil is not required to achieve the same order of accuracy. 
In Figs.~\ref{fig:results_circle}--\ref{fig:results_hourglass_mixed}, it is shown that the problems solved using the $\mathghost$-FEM show worse conditioning compared to those solved with the Coco-Russo scheme. In those tests, the snapping exponent is $\alpha = 2$. The improvement becomes evident in Fig.~\ref{fig:cA_snapping}, where different snapping exponent values are explored for the $\mathghost$-FEM. In that case, the second-order accuracy of the numerical solution is still maintained, particularly for sufficiently large values of $N$, see Fig.~\ref{fig:error_snapping}.

\section{Conclusions}
This study provides a comprehensive comparison between the Coco-Russo scheme, a finite-difference based approach, and the $\mathghost$-FEM, a finite-element based strategy, for solving the Poisson equation in arbitrary domains.
Various geometries and boundary conditions are investigated to evaluate the performance of the two different numerical schemes. The Coco-Russo scheme employs a 9-point stencil to interpolate the solution and its first spatial derivatives at the boundary, ensuring second-order accuracy. Overall, while the $\mathghost$-FEM method does not require an enlarged stencil to achieve second-order accuracy, it faces challenges with the conditioning of the linear systems. However, adjusting the snapping exponent improves some of these issues, maintaining accuracy for larger values of $N$. 

Regarding the  convergence analysis of the two numerical schemes,  for the Coco-Russo scheme remains an open problem. Significant progress has been made in one-dimensional cases, where the scheme has been proven to converge with second-order accuracy in the ${L}^\infty$ norm, and partial extensions of these results to two-dimensional cases, under certain conditions. On the other hand, the $\mathghost$-FEM demonstrates more general and robust convergence properties, achieving second-order accuracy in the ${L}^2$ norm for general domains and boundary conditions. 

\section*{Acknowledgments}
This work has been supported by the Spoke 1 Future HPC \& Big Data of the Italian Research Center on High-Performance Computing, Big Data and Quantum Computing (ICSC) funded by MUR Missione 4 Componente 2 Investimento 1.4: Potenziamento strutture di ricerca e creazione di “campioni nazionali di R \&S (M4C2-19)” - Next Generation EU (NGEU).

The work of A.C.~has been supported from Italian Ministerial grant PRIN 2022 “Efficient numerical schemes and optimal control methods for time-dependent partial differential equations”, No. 2022N9BM3N - Finanziato dall’Unione europea - Next Generation EU – CUP: E53D23005830006; and from the Italian Ministerial grant PRIN 2022 PNRR “FIN4GEO: Forward and Inverse Numerical Modeling of hydrothermal systems in volcanic regions with application to geothermal energy exploitation”, No. P2022BNB97 - Finanziato dall’Unione europea - Next Generation EU – CUP: E53D23017960001.

C.A. and A.C. are members of the Gruppo Nazionale Calcolo Scientifico-Istituto Nazionale di Alta Matematica (GNCS-INdAM).

We acknowledge the CINECA award under the ISCRA initiative (ISCRA C project BCMG, code HP10C7YOPZ), for the availability of high-performance computing resources and support.


%
%



\bibliographystyle{plain}
\bibliography{bibliography.bib}

\end{document}